\documentclass[11pt, twoside]{article}

\pdfoutput=1

\usepackage[all]{xy}

\usepackage{graphicx} 	
\usepackage{subfigure}	
\usepackage{float}		

\usepackage{braket}
\usepackage{amsmath}
\usepackage{amsthm}
\usepackage{amsfonts}
\usepackage{amssymb}
\usepackage{wasysym}

\usepackage[cal=boondoxo]{mathalpha}
\usepackage[toc,page]{appendix}

\usepackage{relsize}
\usepackage{multicol}
\usepackage{multirow}

\usepackage{color}		
\usepackage{tcolorbox}	

\usepackage{latexsym}
\usepackage{mathrsfs}
\usepackage[psamsfonts]{eucal}

\usepackage[margin=1in,paperwidth=8.5in,paperheight=11in]{geometry}

\usepackage[anticlockwise]{rotating}

\usepackage{cancel} 
\usepackage{pdfpages} 
\usepackage{makecell} 
\usepackage{soul} 
\usepackage{enumerate} 
\usepackage{amsthm} 

	\theoremstyle{theorem}
	\newtheorem{theorem}{Theorem}[section] 
	\theoremstyle{theorem}
	\newtheorem{proposition}[theorem]{Proposition}
	
	\theoremstyle{theorem}
	\newtheorem{lemma}[theorem]{Lemma}
	
	\theoremstyle{theorem}
	\newtheorem{corollary}[theorem]{Corollary}
	
	\theoremstyle{theorem}

	\theoremstyle{definition}
	\newtheorem{definition}[theorem]{Definition}
	
	\theoremstyle{definition}

	\theoremstyle{definition}

	\theoremstyle{definition}

	\theoremstyle{definition}

	\theoremstyle{remark}
	\newtheorem{remark}[theorem]{Remark}
\usepackage{imakeidx} 
	\makeindex[intoc]
\usepackage{algorithm} 
\usepackage{algpseudocode} 
\usepackage{array}
\usepackage{mathrsfs} 

\usepackage{tikz} 
\usepackage{tikz-cd} 

\usepackage{mathtools} 
\usepackage{pgfgantt} 
\usepackage{xcolor,colortbl} 
	\definecolor{Grey}{gray}{0.9}
	\definecolor{White}{gray}{1}

\usepackage[T1]{fontenc}
\usepackage{listings} 
\usepackage{wrapfig}
\usepackage{setspace}

\usepackage{ytableau}
\usepackage{fancyhdr}
\usepackage{hyperref}
\hypersetup{pdfborder={0 0 0}}
\usepackage{footnotebackref}

\newcommand{\Hom}{{\rm Hom}}

\newcommand{\codim}{\rm codim}

\newcommand{\Sym}[1]{{\rm Sym}^{#1}}

\newcommand{\Spec}{{\rm Spec}}

\newcommand{\End}{\mbox{End}}

\DeclareMathOperator*{\medotimes}{\raisebox{0.25ex}{\scalebox{0.8}{$\bigotimes$}}}
\DeclareMathOperator*{\medwedge}{\raisebox{0.25ex}{\scalebox{0.8}{$\bigwedge$}}}
\DeclareMathOperator*{\medoplus}{\raisebox{0.25ex}{\scalebox{0.8}{$\bigoplus$}}}

\newcommand{\BBA}{{\mathbb A}}
\newcommand{\BBC}{{\mathbb C}}
\newcommand{\BBR}{{\mathbb R}}

\newcommand{\BBZ}{{\mathbb Z}}

\newcommand{\BBP}{{\mathbb P}}

\newcommand{\BBS}{{\mathbb S}}

\usepackage[backend=biber,style=alphabetic,doi=false,isbn=false,url=false,sorting=nyt]{biblatex}
\renewbibmacro*{finentry}{%
	\setunit{\addperiod\newline \textbf{Annotation:} }%
	\printfield{annotation}%
	\finentry}

\title{{ \Large \bf{HOMOLOGICAL PROJECTIVE DUALITY FOR THE PL{\" U}CKER EMBEDDING OF THE GRASSMANNIAN}}}

\author{BRADLEY DOYLE}

\date{}

\begin{document}


\maketitle

\begin{abstract}
	We describe the Kuznetsov component of the Pl{\" u}cker embedding of the Grassmannian as a category of matrix factorizations on an noncommutative crepant resolution (NCCR) of the affine cone of the Grassmannian. We also extend this to a full homological projective dual (HPD) statement for the Pl{\" u}cker embedding. \par
	The first part is finding and describing the NCCR, which is also of independent interest. We extend results of {\v S}penko and Van den Bergh to prove the existence of an NCCR for the affine cone of the Grassmannian. We then relate this NCCR to a categorical resolution of Kuznetsov. Deforming these categories to categories of matrix factorizations we find the connection to the Kuznetsov component of the Grassmannian via Kn{\"o}rrer periodicity. In the process we prove a derived equivalence between two different NCCR's; this shows Hori duality for the group $SL$. Finally we put this all into the HPD framework.

\end{abstract}

\newpage

\tableofcontents

\section{Introduction}\label{sec:introduction}

The main aim of this paper is finding a homological projective dual (HPD) for the Grassmannian embedding into projective space using the Pl{\" u}cker embedding.
$$ Gr(k,V) \hookrightarrow \BBP\left(\wedge^{k} V \right) $$
The heart of this result is a description of the "interesting"/Kuznetsov component, $\mathcal{K}_H$, of the derived category of a hyperplane section of $Gr(k,V)$. We show it is equivalent to matrix factorizations on a noncommutative crepant resolution (NCCR) of the affine cone of the Grassmannian. \par
First we find an NCCR, $\Lambda$, for the affine cone of the Grassmannian.
\begin{theorem}
	The algebra $\Lambda = \End \left((U \medotimes \Sym{}(\Hom(S,V)^*))^{SL(S)}\right)$ is an NCCR for the affine cone of the Grassmannian.
\end{theorem}
Here $U$ is a certain $SL(S)$ representation, see Theorem \ref{thm:FGDandNCCR} for the exact statement. This is a result both of independent interest and is also needed to describe the HPD. Next we relate this NCCR to a certain subcategory of $D^b(\mathcal{O}_{Gr(k,V)}(-1))$ and in the process we also show a derived equivalence between two different NCCR's. This derived equivalence is a special case of Hori duality. Building on these results, we relate $\mathcal{K}_H$ to a category of matrix factorizations,  $D^b(\mathcal{O}_{Gr(k,V)}(-1), W)$, using Kn{\" o}rrer periodicity.
The result we get is
$$\mathcal{K}_H \simeq D^b(\Lambda, W)$$
where $W$ is a superpotential related to the hyperplane $H$. We finally extend this to a full HPD statement.
\begin{theorem}
	The HPD for the Pl{\" u}cker embedding of the Grassmannian is the category of matrix factorizations on an NCCR for a fibre bundle over $\BBP(\wedge^k V^*)$. Over a single fibre, this NCCR is isomorphic to $\Lambda$.
\end{theorem}
See Theorem \ref{thm:finalequiv} for a precise statement. In the rest of the introduction we will give an extended summary, providing motivation and a sketch of the results and techniques used as well as references.

\subsection{Non-commutative crepant resolutions and categorical resolutions}\label{sec:nccrintro}
We begin with a summary of non-commutative resolutions. These are often required in HPD, as the varieties arising are not smooth, but are also an independent subject on their own. \par
The idea of trying to resolve a singularity by replacing a space which is not smooth with another space birational to the original space and smooth has been very useful in the study of schemes. This is a well known and very commonly used tool in algebraic geometry. As resolutions of singularities are not unique we often ask for extra conditions, one of those conditions is crepancy, which means the canonical bundle pulls back to the canonical bundle. This type of resolution is geometric, recently the idea of taking a different type of resolution has arisen, using non-commutative rings or categories instead. This is no longer a strictly geometric idea, however it still has connections to geometry and even to classical resolutions of singularities in certain cases. \par
The first mentions of non-commutative resolutions were by Bondal and Orlov \cite{Bondal02derivedcategories} and Van den Bergh \cite{vandenbergh2004}. Bondal and Orlov defined a categorical desingularization of $X$ as a category with finite homological dimension which had a quotient equal to $D^b(X)$. Van den Bergh studied endomorphism algebras of modules over a singular ring. These algebras can be thought of as a noncommutative resolution. Both of these papers originally saw these resolutions as something that could be useful to show other results. Van den Bergh only defines an non commutative (crepant) resolution in the appendix of \cite{vandenbergh2004}, it had been used as an intermediary step in proving two derived categories were equivalent. In \cite{van2004non}, Van den Bergh studies non commutative (crepant) resolutions, in particular their existence in dimension $3$. Looking at derived categories and equivalences was also the main content of the Bondal and Orlov paper.\par
Later, {\v S}penko and Van den Bergh, \cite{2015arXiv150205240S}, gave a framework for proving the existence of NCCR's for quotients of smooth varieties by a reductive group, $G$. They prove the existence of NCCR's under specific conditions, and it is their definition that we use.
\begin{definition}[{\cite[Def. 1.1.1]{2015arXiv150205240S}}]\label{def:nccr}
	Let $\Spec \;S$ be a affine scheme, such that $S$ is a normal noetherian domain. \\
	A \emph{non-commutative resolution} of $S$ is an algebra over $S$ which has finite global dimension and is of the form $\Lambda = \End_S(M)$, for $M$ a non-zero, finitely generated, reflexive $S$-module. It is called \emph{crepant} if $S$ is Gorenstein and $\Lambda$ is a maximal Cohen-Macaulay $S$-module.\\
	We abbreviate with NC(C)R. 
\end{definition}
See Van den Bergh \cite[Def 4.1]{van2004non} for the original definition and Leuschke \cite{2011arXiv1103.5380L} for motivation and more details about this definition. \par
Kuznetsov, \cite{2006math......9240K}, defined a type of categorical resolution and proved an existence result, again under specific conditions. Later, Kuznetsov with Lunts, \cite{2012arXiv1212.6170K}, proved the existence of a categorical resolution for the derived category of any separated scheme with characteristic 0.
\begin{definition}[{\cite[Def. 3.2]{2006math......9240K}}]\label{def:categoricalresolution}
	A \emph{categorical resolution} of a noetherian scheme $S$ is a regular triangulated category $T$ and maps
	$$\pi_{*} : T \rightarrow D^b(S) \qquad \qquad \pi^*: D^b(S) \rightarrow T$$
	such that $\pi^*$ is the left adjoint to $\pi_*$ on $D^{perf}(S)$ and the natural map $id_{D^{perf}(S)} \rightarrow \pi_* \pi^*$ is an isomorphism. \\
	It is \emph{weakly crepant} if $\pi^*$ is also the right adjoint to $\pi_*$ on $D^{perf}(S)$. If, in addition,  $T$ has a module structure over $S$ and the relative Serre functor is the identity, it is \emph{strongly crepant}. 
\end{definition}
Note that if $S$ has rational singularities, then any geometric resolution $\pi:\tilde{S} \rightarrow S$ gives a categorical resolution $D^b(\tilde{S})$ and in this setting, weakly crepant, strongly crepant and $\pi$ being a crepant geometric resolution are all equivalent. See \cite{abuaf2016categorical} for more details on crepancy for categorical resolutions. \par
An NC(C)R gives a categorical (strongly crepant) resolution and if the categorical resolution is generated by a tilting bundle, the reverse is also true. In fact categorical resolutions are more general and are probably close to the correct concept of resolutions for categories, but NC(C)R's are more concrete and easy to study. Also both \cite{2006math......9240K} and \cite{2012arXiv1212.6170K} use existing geometric resolutions to find categorical resolutions, whereas \cite{2015arXiv150205240S} does not. For more in depth motivation and discussion see the survey paper by Leuschke \cite{2011arXiv1103.5380L}.\par
{\v S}penko and Van den Bergh take motivation from the case of quotients by finite groups, which can be briefly summarised as follows: \\
Let $G$ be a finite group acting on a regular ring $R$. We can define the twisted group ring, $G\#R$,  modules over this ring are exactly the $G$-equivariant $R$-modules. This ring has finite global dimension equal to the global dimension of $R$ and in fact it is an NC(C)R of $R/G$. There is an alternative way of describing $G\#R$. Let $M$ be the module of covariants associated to the sum of all the irreducible representations, then $G\#R \cong \End_{R^G}(M)$. This description motivates the infinite group case.\par
In the case of an infinite group, you can not take all the irreducible representations, so {\v S}penko and Van den Bergh take a finite collection of irreducible representations. Then they study potential NC(C)R's of the form
$$\Lambda(M) = \End_{R^G}(M)$$
where $M$ is a finite sum of modules of covariants. They then construct resolutions related to general modules of covariants and using these resolutions, prove conditions for $G$ actions and choices of $M$, under which $\Lambda(M)$ has finite global dimension.\par
Earlier,  Van den Bergh \cite{van1991cohen}, proved results about the Cohen-Macaulayness of modules of covariants, they also show under what conditions these results give Cohen-Macaulayness of $\Lambda(M)$.\par
\vspace{1cm}
In Section \ref{sec:algebraicnccr} we will find an NCCR for a family of examples, the affine cone of the Grassmannian. Let $S$ and $V$ be $k$ and $n$ dimensional vector spaces over $\BBC$ and let $X = \Hom (S, V)$, where $n - 1 > k > 1$. Consider the standard action of $G = SL(S)$ on $X$, then $X/G$ is the affine cone of the Grassmannian, $Gr(S,V)$. It has one singular point, the origin, and has a natural geometric resolution $\pi: \mathcal{O}(-1)_{Gr(S,V)} \rightarrow X/G$ . We want to find an NCCR of $X/G$ using the framework of \cite{2015arXiv150205240S}. We require that $n,k$ are coprime.\footnote{See subsection \ref{sec:HPDintro} for more details about this condition.}\par
We can not just apply the main result of \cite{2015arXiv150205240S} as our example is not quasi-symmetric, so we will need to strengthen results about Cohen-Macaulayness and finite global dimension. First of all we deal with Cohen-Macaulayness, going back to \cite{van1991cohen}, we analyse the spectral sequence provided there and show which modules of covariants are Cohen-Macaulay. Second, we explain how a result of Fonarev \cite{fonarev2013minimal} can be easily adapted to show finite global dimension.\par
Putting these results together in Theorem \ref{thm:FGDandNCCR} we get an NCCR
$$\Lambda = \End \left(\medoplus  \left(\BBS^{\alpha}S^{*} \medotimes \Sym{\bullet}X^{*} \right)^G\right)$$
where $\BBS^{\alpha}S^{*}$ is a Schur functor and the sum is over all Young diagrams that fit inside a triangle of length $n-k$ and height $k$, see Section \ref{sec:schurandyoung}. \par
Consider another vector space, $Q$, of dimension $n-k$. We get an isomorphic singularity, $\Hom(V, Q)/ SL(Q)$, and a very similar NCCR, $\Lambda'$. In general however, the algebras $\Lambda$ and $\Lambda'$ are not isomorphic. We prove instead that these algebras are derived equivalent, see Theorem \ref{thm:allNCCRsareequiv}. \par
This gives a mathematical proof of a version of Hori duality for the group $SL$, at least for $n$ and $k$ co-prime. See \cite{2007JHEP...05..079H} and \cite{2011arXiv1104.2853H}. The brief summary is that for theoretical physics reasons there should be a duality on a gauge theory level between $SU(k)$ and $SU(n-k)$ as $Gr(k,n) = Gr(n, n-k)$. Mathematically it is thought that $D^b(\Lambda)$ and $D^b(\Lambda')$ should be models for the category of B-branes, which are special types of boundary conditions for the gauge theory. Therefore this equivalence is predicted by Hori duality, assuming they are actually mathematical models for the B-branes. See also \cite{2016arXiv160904045V} where Rennemo and Segal prove Hori duality for the Orthogonal group.\par
This derived equivalence between $\Lambda$ and $\Lambda'$ will follow as a corollary of the proof that both NCCR's are derived equivalent to a categorical resolution. This categorical resolution is found by applying the work of Kuznetsov \cite{2006math......9240K} to the geometric resolution, $\pi: \mathcal{O}(-1)_{Gr(S,V)} \rightarrow X/G$. To prove this equivalence of categories we use the stack \\ $\left[ \Hom(S,V) \times \wedge^k S/ GL(S)\right]$, both $\left[X/G\right]$ and $\mathcal{O}_{Gr(S,V)}(-1)$ are open sets inside this stack, and we will use a "window" subcategory style argument to relate the NCCR to the categorical resolution. Section \ref{sec:equivNCCR} deals with this. This equivalence also forms the heart of the equivalence $\mathcal{C}_H \simeq D^b(\Lambda, W)$ which itself leads to our HPD result. See the rest of the introduction for more details.\par
We finally note that even though Kuznetsov's categorical resolutions and {\v S}penko-Van den Bergh's NCCR are different (but related) ideas, the proof for finite global dimension of the NCCR is effectively the same as the proof that gives us generation in the Lefschetz decomposition. It would not be surprising if this turned out to be true in other situations as well.

\subsection{Hyperplane sections of the Grassmannian}\label{sec:hyperplanesections}
Before considering the Grassmannian, we look at the simpler case of $\BBP^{n}$. Beilinson, \cite{beilinson1978coherent}, tells us that
$$D^b(\BBP^n) = \langle \mathcal{O}, \mathcal{O}(1), \dots , \mathcal{O}(n) \rangle.$$
Now consider a degree $d \leq n+1$ hypersurface $X_d \subset \BBP^n$. It is cut out by a section of $\mathcal{O}(d)$ and we have the exact sequence
$$0 \rightarrow \mathcal{O}_{\BBP^n}(-d) \rightarrow \mathcal{O}_{\BBP^n} \rightarrow \mathcal{O}_{X_d} \rightarrow 0.$$
One can use this to prove that the restriction of $\mathcal{O}(d), \dots, \mathcal{O}(n)$ to $X_d$ is an exceptional collection. This tells us that
$$ D^b(X_d) = \langle \mathcal{K}_d, \mathcal{O}(d), \dots , \mathcal{O}(n) \rangle $$
where $\mathcal{K}_d$ is the right orthogonal to $ \langle \mathcal{O}(d), \dots , \mathcal{O}(n) \rangle$ in $D^b(X_d)$, it is called the "interesting" or Kuznetsov component. From a derived category point of view, to understand $X_d$ it is sufficient to understand $\mathcal{K}_d$. \par
If $d=1$, $X_d \cong \BBP^{n-1}$ and $\mathcal{K}_1 \simeq 0$. If $d=2$, then $X_d$ is a quadric and $\mathcal{K}_2$ is generated by 2 completely orthogonal objects if n is even, and 1 object if $n$ is odd. For $3 \leq d \leq n+1$ the description is more complicated. Fix $2 \leq d \leq n+1$, then we have 
\begin{center}
	\[\begin{tikzcd}
	& {\mathcal{O}(-d)} \\
	{X_d} & {\BBP^n}
	\arrow["p"', from=1-2, to=2-2]
	\arrow["i", hook, from=2-1, to=2-2]
	\end{tikzcd}\]
\end{center}
and $D^b(X_d) \simeq D^b(\mathcal{O}(-d), W)$ via Kn{\"o}rrer periodicity. Here $D^b(\mathcal{O}(-d), W)$ is a category of matrix factorizations and $W: \mathcal{O}(-d) \rightarrow \BBC$ is the superpotential, it is defined using the section whose zero locus is $X_d$. See Section \ref{sec:background} for more information on the various mathematics used in the introduction. Using Kn{\"o}rrer periodicity we can find a subcategory $\mathcal{D}_d \subset D^b(\mathcal{O}(-d), W)$ that is equivalent to $\mathcal{K}_d$.\par
Next, note that $\mathcal{O}(-d)$ and $\left[ \BBA^{n+1}/\BBZ_d \right]$ are both resolutions of $\BBA^{n+1}/\BBZ_d$, where $\BBZ_d$ acts coordinate-wise by multiplication with a $d^{th}$ root of unity. We can use VGIT to relate them, in general we get an inclusion\footnote{In the special case $d=n+1$ this is a flop and we get an equivalence of categories.}
$$D^b\left(\left[\BBA^{n+1}/\BBZ_d\right], W\right) \hookrightarrow D^b(\mathcal{O}(-d), W).$$
It turns out that the image of this inclusion is equivalent to $\mathcal{K}_d$ and therefore we have
$$\mathcal{K}_d \simeq D^b\left(\left[\BBA^{n+1}/\BBZ_d\right], W\right).$$
From Section \ref{sec:nccrintro} we can view $D^b(\left[\BBA^{n+1}/\BBZ_d\right])$ as an NC(C)R of the singularity $\BBA^{n+1}/\BBZ_d$. In particular, $D^b(\left[\BBA^{n+1}/\BBZ_d\right])$ is just the derived category of the twisted group ring, and we have an action by a finite group,  so this twisted group ring is an NC(C)R. Therefore we can say that the Kuznetsov component of $X_d$ is matrix factorizations on a specific NC(C)R. It is crepant only when $d$ divides $n+1$. See \cite{abuaf2016categorical} for details on the NC(C)R part and \cite{2013arXiv1306.3957B} for more details on the rest, that paper is all written from the homological projective duality point of view, looking at the Veronese map for projective space, we will introduce homological projective duality in the next section. It is also these ideas that we will generalise to the Grassmannian, again see the next section for more details.\par
The Kuznetsov component exists in greater generality than just for projective space. Let $X$ be a smooth projective variety with a line bundle, $\mathcal{L}$, such that
$$D^b(X) = \langle \mathcal{A}_0, \mathcal{A}_1(\mathcal{L}), \dots , \mathcal{A}_n(\mathcal{L}^{\otimes n}) \rangle$$
where $\mathcal{A}_n \subseteq \mathcal{A}_{n-1} \subseteq \cdots \subseteq \mathcal{A}_0$ are admissible triangulated categories. This is called a {\it Lefschetz decomposition}, if in addition $\mathcal{A}_0 = \mathcal{A}_1 = \dots = \mathcal{A}_n$ we call it a \emph{rectangular Lefschetz decomposition}. Then let $Y \subset X$ be a hypersurface cut out by a regular section of $\mathcal{L}$, we get
$$D^b(Y) = \langle \mathcal{K}, \mathcal{A}_1(\mathcal{L}), \dots , \mathcal{A}_n(\mathcal{L}^{\otimes n}) \rangle $$
and again, we wish to understand the Kuznetsov component, $\mathcal{K}$.\par
Of particular interest is the case where the $\mathcal{A}_i$ are generated by an exceptional collection, i.e. $X$ has a full exceptional collection. If this is the case, then if $\mathcal{K}$ is trivial or itself generated by exceptional objects we get that $Y$ has a full exceptional collection as well.\par
We will look at the Grassmannian. Kapranov showed that we have a full exceptional collection and in the special case, $Gr(2,n)$, Kuznetsov \cite{kuznetsov2005expceptionalgrassmannian}, found a minimal\footnote{One can put an ordering on Lefschetz decompositions, given by the inclusion relation on $\mathcal{A}_0$.} Lefschetz collection. Fonarev, \cite{2011arXiv1108.2292F}, then generalised this to $Gr(k,n)$, assuming that $k$ and $n$ are coprime, he found a rectangular Lefschetz decomposition where $\mathcal{A}_0$ is generated by a collection of vector bundles associated to representations of $GL(k)$. \par
In \cite{kuznetsov2005expceptionalgrassmannian}, Kuznetsov showed that for $Gr(2,2n)$, $\mathcal{K}$ is trivial, and his method extends to $Gr(2,2n+1)$, again in this case $\mathcal{K}$ is trivial.\footnote{Use his method to relate hyperplane sections of $Gr(2,2m)$ and $Gr(2,2m+1)$, then use the result for $Gr(2,2m)$ to conclude the result for $Gr(2,2m+1)$.} \par
We will show that for $Gr(k,n)$,
$$\mathcal{K} \simeq D^b(\Lambda, W)$$
where $\Lambda$ is the NCCR from earlier on in the introduction and $W$ is a superpotential defined later. For more details see the next section, as this is best seen as the heart of the HPD result.
\subsection{Homological Projective Duality}\label{sec:HPDintro}
Instead of just considering a single hyperplane section, and trying to understand the Kuznetsov component, $\mathcal{K}$, homological projective duality (HPD) seeks to understand the categories $\mathcal{K}_H$ arising from all possible choices of linear hyperplanes, $H$, as well as higher codimension linear subspaces. In other words it is generalizing $\mathcal{K}$ to a family over the space of hyperplanes, which is just the dual projective space. See \cite{2005math......7292K} for Kuznetsov's original work on HPD and \cite{2018arXiv180400132P} for Perry's generalization of HPD to noncommutative algebraic geometry. HPD can be thought of as an adaptation of projective duality to categories.\par
The basic idea is to embed a projective variety, $X$, into projective space 
$$X \rightarrow \BBP V^*$$
and then find another variety, $Y$, defined over the dual projective space, $\BBP V$, such that $D^b(Y)$ and $D^b(X)$ are related in a specific way. In particular, assume that $D^b(X) = \langle \mathcal{A}, \dots, \mathcal{A}(n) \rangle$ where $\mathcal{O}_X(1)$ is the pullback of $\mathcal{O}_{\BBP V^*}(1)$. Then take a point in $\BBP V$, it corresponds to a hyperplane, $H$, and by intersecting $X$ with $H$ we get $D^b(X_H)$, which decomposes as follows
$$D^b(X_H) = \langle \mathcal{K}_H, \mathcal{A}(1), \dots, \mathcal{A}(n) \rangle.$$
We can also consider the fibre of $Y$ at the point $H$, denote it $Y_H$. We want $D^b(Y_H) \cong \mathcal{K}_H$.\par
Kuznetsov, \cite[Definiton 6.1]{2005math......7292K} defines a HPD as a variety $Y \rightarrow \BBP V$ such that the universal hyperplane section $\mathcal{X}_1 \subset X \times \BBP V$ has a specific semi-orthogonal decomposition containing $\phi \left(D^b(Y)\right)$ where $\phi$ is a fully faithful Fourier-Mukai map. (This is a formalization of the above discussion). Then Kuznetsov proves the following collection of equivalences holds.\par
Let $V = H^0(X, \mathcal{O}(1))$ and let $L \subset V$ be any linear subspace, set $L^{\perp} \subset V^*$ to be the annihilator of $L$. Then $X \hookrightarrow \BBP V^*$ and we set $X_{L^{\perp}}$ to be the intersection of $X$ and $L^{\perp}$. We get $Y_L \rightarrow \BBP L$ and an inclusion of subcategories between $D^b(X_{L^{\perp}})$ and $D^b(Y_L)$, which one includes into the other depends on $\dim L$. 
\begin{theorem}[Fundamental Theorem of Homological Projective Duality]\label{thm:HPDtheorem}
	Let $X \hookrightarrow \BBP V^*$ and $\mathcal{O}_X(1)$ be the pullback of $\mathcal{O}_{\BBP V^*}(1)$. Assume $D^b(X) = \langle \mathcal{A}, \mathcal{A}(1), \dots , \mathcal{A}(n) \rangle$. Let $L \subset V$ be a linear subset with $\dim L = l$, set $L^{\perp} \subset V^*$ to be the annihilator of $L$. Let $Y \rightarrow \BBP V$ be an HPD for $X$, and let $\mathcal{L}$ be the pullback of $\mathcal{O}_{\BBP V}(1)$. Set $X_{L^\perp}$ and $Y_L$ to be the base changes of $X$ and $Y$. We have the following equivalences
	\begin{itemize}
		\item If $l < n + 1, \qquad D^b\left( X_{L^{\perp}}\right) \simeq \langle D^b(Y_L), \mathcal{A}(l), \dots , \mathcal{A}(n) \rangle$ \\
		\item If $l = n + 1, \qquad D^b\left( X_{L^{\perp}}\right) \simeq D^b(Y_L)$ \\
		\item If $l > n + 1, \qquad D^b(Y_L) \simeq \langle \mathcal{A}((- l + 1)\mathcal{L}), \dots , \mathcal{A}(- n\mathcal{L}),  D^b\left(X_{L^{\perp}}\right)\rangle$
	\end{itemize}
	whenever both $X_{L^\perp}$ and $Y_L$ have the expected dimension. 
\end{theorem}
\begin{remark}
	The $\mathcal{A}$ that appears in the decomposition of $D^b(Y_L)$ is more accurately the image of $\mathcal{A}$ along a fully faithful functor. In specific, $\mathcal{A}$ embeds into $D^b(Y)$ and can then be projected into $D^b(Y_L)$.
\end{remark}
\begin{remark}\label{rmk:dim1case}
	If $\dim L = 1$, we get out $\mathcal{K}_H \simeq D^b(Y_L)$. Also, if $L = V$, then we get that $D^b(Y)$ has a rectangular Lefschetz decomposition, this is often stated as another part of the fundamental theorem. Finally, this is actually a duality, one can prove that $X$ is HP dual to $Y$.
\end{remark}
In practice it turns out that it is not normally possible to find a \emph{variety} with this property, so we have to expand to less geometric objects. If one lets $Y$ be a triangulated category linear over the dual projective space, one can always find a HPD (Kuznetsov \cite{2005math......7292K}, Perry \cite{2018arXiv180400132P}). However for this to be interesting, we seek to find some "geometric" description of $Y$. In some cases it is locally a variety, but with a global Brauer twist, sometimes a NCCR of a variety, maybe a matrix factorization category associated to a variety or NCCR. There is also a similar statement for the more general case of an arbitrary Lefschetz decomposition.\par
We however will not use Kuznetsov's definition, but will instead define a HPD  as any triangulated category linear over $\BBP V$ that Theorem \ref{thm:HPDtheorem} is true for. This partially fits in with \cite{2016arXiv160904045V}, \cite{2017arXiv170501437V} and \cite{2013arXiv1306.3957B}. In these papers, there are different but "equivalent" definitions, you can find discussions about their choices and the relationships between them in those papers.\par
Our approach is closest to the framework studied by Rennemo, \cite{2017arXiv170501437V}. He works more generally with subcategories $\mathcal{D} \subset D^b([Z/G])$ which have a Lefschetz decomposition and replaces $D^b(X_{L^\perp})$ and $D^b(Y_L)$ with triangulated categories, $\mathcal{D}_{L^{\perp}}$ and $\mathcal{D}^*_L$. In this setting we do not have to worry about having the expected dimension as $\mathcal{D}_{L^\perp}$ and $\mathcal{D}^*_L$ should be thought of as the correct categorical base changes.\par
In more detail, let $\mathcal{X}_1 \subset \mathcal{Z} = \mathcal{O}_{X \times \BBP V}(-1,-1)$ be the universal hyperplane section. Then $D^b(\mathcal{Z}, W) \simeq D^b(\mathcal{X}_1)$ via Kn{\"o}rrer periodicity. We have an inclusion $\iota_v: X \rightarrow \mathcal{Z}$ for any nonzero  $v \in V$ and then \cite{2017arXiv170501437V} defines the HPD as a certain subcategory, $\mathcal{D}^{*} \subset D^b(\mathcal{Z}, W)$. An object, $\mathcal{E}$, is an element of $\mathcal{D}^{*}$ if $\iota^*_v \mathcal{E} \in \mathcal{A}_0$. Our HPD will be defined similarly and as part of our proof we will show an equivalence to Rennemo's HPD, proving that we do have an HPD.\par
We also use ideas and methods from Segal, \cite{2011CMaPh.304..411S}, Ballard-Deliu-Favero-Isik-Katzarkov, \cite{2013arXiv1306.3957B} as well as others.\par
In more detail, we consider the Pl{\" u}cker embedding of the Grassmannian. Let $V$ be a vector space of dimension $n$, pick $k$ such that $1<k<n-1$ and $n$ and $k$ are coprime. Let $X = Gr(k,V)$, we have
\[\begin{tikzcd}
{X} & {\BBP(\wedge^k V)}
\arrow[from=1-1, to=1-2, hook]
\end{tikzcd}\]  
and the Kuznetsov - Fonarev exceptional collection
$$D^b(X) = \langle \mathcal{A}, \dots ,\mathcal{A}(n-1) \rangle$$
where $\mathcal{A}$ is the set of vector bundles associated to a specific collection of $GL(k)$ representations, see Section \ref{sec:schurandyoung} for more details.\par
Looking at the case of a single hyperplane section again, we have $X_H$, i.e. $X$ and a section $w$ of $\mathcal{O}_X(1)$ whose zero locus is $X_H$. Kn{\" o}rrer periodicity tells us that $D^b(X_H) \cong D^b(\mathcal{O}_X(-1), W)$, where the superpotential is $W = pw$ for $p$ the fibre coordinate. See Section \ref{sec:stackandmf} for more details on matrix factorizations and Kn{\" o}rrer periodicity. \par
As we are interested in $\mathcal{K}_H \subset D^b(X_H)$ we want to find the corresponding subcategory of $D^b(\mathcal{O}_X(-1), W).$ Forgetting $W$, inside $D^b(\mathcal{O}_X(-1))$ there is a subcategory, $\mathcal{B}$, a categorical resolution of the affine cone of the Grassmannian, see Section \ref{sec:nccrintro} above and Section \ref{sec:equivNCCR}. As part of the proof that the two algebraic NCCR's are equivalent we prove that they are both equivalent to $\mathcal{B}$. As the definition of $\mathcal{B}$ is closely related to the definition of the HPD that Rennemo gives, this suggests that the HPD should be closely related to the algebraic NCCR with a superpotential added on. We split the general HPD result into three steps, each of which we deal with separately in Section \ref{sec:HPD}.
\begin{itemize}
	\item Find the subcategory of $D^b(\mathcal{O}_X(-1), W)$ that is equivalent to $\mathcal{K}_H$.
	\item Extend the equivalence between the algebraic and categorical resolution to include the superpotential $W$. This will prove $\mathcal{K}_H \simeq D^b(\Lambda, W)$.
	\item Extend the above equivalence to a HPD statement. (i.e. all the statements of equivalences for $L \subset \wedge^k V^*$.)
\end{itemize}
The first point is straightforward in the case of a single hyperplane, the general version over $\BBP L \subset \BBP \left( \wedge^k V^*\right)$ will follow from work of Rennemo \cite{2017arXiv170501437V} and the third point is more of an important technical statement, but conceptually it follows from the equivalence between the algebraic and categorical resolutions. See Propositions \ref{prop:linebundleequivfibre} and \ref{prop:nccrequivfibre} for the results for a single fibre, Theorem \ref{thm:linebundleenccrquivpotential} for the matrix factorization version and finally Theorems \ref{thm:secondequiv} and \ref{thm:finalequiv} for the main result.
\begin{remark}
	We can alternatively describe the HPD as follows: Geometrically we have a family of singularities over $\BBP (\wedge^k V^*)$, each fibre is isomorphic to $\Hom(S,V)/SL(S)$. Then forgetting the superpotential $W$, the HPD is the NCCR of this singular space. Over any point of $\BBP (\wedge^k V^*)$, this NCCR is exactly the NCCR, $\Lambda$ from Section \ref{sec:nccrintro}, however the potential varies over each point. In fact, the point effectively defines the superpotential.
\end{remark}
\begin{remark}
	Our geometric understanding of the HPD still is not very strong. Both individual parts that make it up have good geometric and algebraic understanding. NCCR's are a natural concept closely related to geometry, and a subcategory of a derived category of a quotient stack generated by vector bundles is also not that difficult to understand. The other part of the HPD, matrix factorizations is also well studied, and is a model for some ideas from theoretical physics. In some situations these objects are actually equivalent to schemes. (i.e, Kn{\" o}rrer periodicity.) However when we combine these ideas together and get matrix factorizations on a NCCR, we do not have much geometric information any more and there currently are not many results about these objects.\par
	We do not expect a much better geometric description in general, but for specific $k$ and $n$ there may be one. In particular, Kuznetsov looked at HPD for $Gr(2,n)$ for $n \leq 7$ and Rennemo and Segal, \cite{2016arXiv160904045V}, looked at HPD for the Pfaffian, of which a special case is $Gr(2,n)$. In these cases the HPD is supported over a subvariety of dual projective space as $\mathcal{K}$ is generically zero, see the start of this section. This subvariety turns out to be another Pfaffian. Both sides also require a non-commutative resolution.
\end{remark}

In all three sections of the introduction we have required that $k$ and $n$ are coprime. The main reason for this is that $Gr(k,n)$ has a rectangular Lefschetz decomposition only in this case. If $k$ and $n$ are not coprime we only have a conjectured minimal Lefschetz decomposition, which is not rectangular, \cite{fonarev2013minimal}.  Without rectangular decompositions we still expect to get categorical, and non-commutative resolutions, but we no longer expect these resolutions to be crepant. See \cite[Section 5]{2016arXiv160904045V} for a brief discussion of the issues in dealing with $Gr(2,2n)$. Some of the calculations in Section \ref{sec:algebraicnccr} rely on the fact that $n$ and $k$ are coprime.

\subsection{$Gr(3,7)$ and other special cases}\label{sec:specialcases}
This final part of the introduction is tangential to the work we do in this paper, but is included here for completeness.\par
Recall the question from Section \ref{sec:hyperplanesections}. "What can we say about $\mathcal{K}$ for a hyperplane section of the Grassmannian?" For arbitrary\footnote{We assume $k < n/2$ as there is symmetry between $k,n$ and $n-k, n$.} $n/2 > k > 2$, coprime, there are a few things we can say. It is a Calabi-Yau (CY) category with CY dimension $k(n-k) + 1 - 2n$, see \cite{2015arXiv150907657K}. If $k \geq 3$ and $n$ is large enough this will be non-zero. We can also calculate the Hodge diamond of the hyperplane section using the Lefschetz hyperplane theorem, see \cite{2019arXiv191203144B} for more details including the formula to make this calculation. Using this we can extract the Hochschild homology of $\mathcal{K}$ as Hochschild homology splits across semi-orthogonal decompositions. In all the low dimensional cases apart from $Gr(2,n)$ and $Gr(3,7)$ we get a non-zero value for $HH_0(\mathcal{K})$. We expect that in all other cases $\mathcal{K} \not \simeq 0$. \par
Apart from $k=2$, there is only one case, $Gr(3,7)$, where we get a negative value for the CY-dimension of $\mathcal{K}$. We also get $HH_0(\mathcal{K})=0$ in this case which suggests $\mathcal{K}$ is trivial. We get CY dimension $0$ for $Gr(3,8)$ and $\dim HH_0(\mathcal{K}) = 2$ which suggests that in this case, $\mathcal{K}$ should be generated by two completely orthogonal objects and therefore the HPD is generically a double cover.\par
There is one case where $\mathcal{K}$ is a CY2 category, $Gr(3,10)$, and this case is studied in detail in \cite{2019arXiv191203144B}, focusing on Hodge theoretic results for the hyperplane section. They prove that for $Gr(3,10)$, the HPD is a non-commutative K3 fibration, but definitely not an actual K3 fibration, in particular $\mathcal{K}$ is not the category of any K3 surface.\par
There is also one other non-coprime case that has been studied, $Gr(3,6)$, in \cite{deliu2011homological}. That thesis provides evidence that suggests the HPD is generically a double cover, similar to the $Gr(3,8)$ case. For higher values of $k$ and $n$ there is not anything interesting that we can currently say.\par
As the numerics suggest that $\mathcal{K} \simeq 0$ for $Gr(3,7)$, we will briefly sketch two possible methods to prove this. This is equivalent to proving that the obvious exceptional collection is in fact full.\par
First, adapting the techniques of \cite{kuznetsov2005expceptionalgrassmannian}. To briefly summarise that paper, Kuznetsov uses induction to find full exceptional collections of hyperplane sections of Grassmannians. Let $Y$ be a hyperplane section of $Gr(2,n)$, and assume a hyperplane section, $Y'$, of a smaller Grassmannian has a full exceptional collection. Then use the Koszul resolution of this smaller Grassmannian to find a full collection on $Y$. This collection turns out to be larger than we want so we need to find sequences relating objects in the collection. Using those sequences, some objects can be removed while still having a full collection. Most of the sequences exists on the Grassmannian, but the key sequence only exists on $Y$. It is closely related to a matrix factorization of the section cutting out $Y$, so restricting to $Y$ gives a complex.\par
Therefore one can look for an analogue of this complex in the $Gr(3,7)$ case. As $Gr(3,7)$ is the only case where we expect $\mathcal{K} = 0$ we are not sure if there are special complexes in this case only, probably for numerical reasons, or if there are general complexes for $Gr(3,n)$, but they only show fullness in the $Gr(3,7)$ case.\par
The second approach uses the geometry of the hyperplane section, $Y \subset Gr(3,7)$. In particular $Y$ has a natural action of $G_2$ as the $3$-form defining $Y$ is exactly the form defining $G_2$. $Y$ is not a homogeneous space for $G_2$, but we can describe it in terms of $G_2$ representations. \\
We have $Gr(3,7) \hookrightarrow \BBP(\wedge^3 V)$ where $\dim V = 7$ and we think of it as the natural $7$ dimensional representation of $G_2$. Then $\wedge^3 V \cong \BBC \medoplus V \medoplus W$ as $G_2$ representations, where $W$ is the unique representation of $G_2$ with dimension 27, it has highest weight $2v_1$ where $V$ has highest weight $v_1$. (In fact $\BBC \medoplus W \cong \Sym{2} V$.) This tells us that $Y \subset \BBP(V \medoplus W)$ and therefore it has two closed orbits, $X_1$ and $X_2$, both of which are isomorphic to $Q^5$, a homogeneous space for $G_2$, and an open orbit of dimension 11, isomorphic to $G_2/SO_3$. \par
In \cite{2018arXiv180305063G} they work with similar objects, called horospherical varieties, unfortunately, this example is not horospherical but a modification of their methods may apply. The idea is to blow up $Y$ at one or both of its closed orbits, then also describe this blowup, $\tilde{Y}$, as a fibre bundle over a $G_2$ homogeneous space $Z$. Then as the closed orbits and $Z$ have known full exceptional collections, if the fibres also have full exceptional collections we would have two descriptions of $D^b(\tilde{Y})$. One from the description as a fibre bundle and one from the fact it is a blowup. Once you have these two descriptions you could attempt to prove that $\mathcal{K}$ was zero using a mutation argument in $\tilde{Y}$. \par
Unfortunately we have not managed to make progress with either strategy.
\begin{remark}
	There is a general conjecture that predicts that any homogeneous space has a full exceptional collection.  Now $Y$ is not a homogeneous space, but still has an action of a reductive group. It would be interesting to see how far this conjecture can be extended. In particular it seems likely to extend to many horospherical varieties, see \cite{2018arXiv180305063G}.
\end{remark}

\subsection{Acknowledgements}\label{sec:acknowledgements}
I would like to thank Ed Segal for supervising me throughout my PhD and always being keen to explain and answer any question I had. I would also like to thank Federico Barbacovi for various conversations about algebraic geometry and Benjamin Aslan for the same about differential geometry. Finally a big thank you to my mom for listening to me attempt to explain my studies for the past 8 years.\par
This work was supported by the Engineering and Physical Sciences Research Council \\
\text{[EP/L015234/1]}. The EPSRC Centre for Doctoral Training in Geometry and Number Theory (The London School of Geometry and Number Theory), University College London.
\subsection{Notation, conventions and assumptions}
All vector spaces, varieties and stacks are over $\BBC$, all functors will be derived unless stated otherwise, we ignore shifts in triangulated categories, in particular the exact triangles are correct up to shifts. \par We use $X$ for a variety and $\mathcal{X}$ for a stack. \par Algebraic groups will be denoted by $G$ or $H$, and $T $ is a maximal torus, $B$ is a Borel subgroup and $P$ is a parabolic subgroup, we assume $T \subset B \subset P$. \par The category of coherent sheaves is coh, the bounded derived category is $D^b$ and triangulated subcategories are denoted by curly letters, $\mathcal{C}$. Triangulated subcategories of matrix factorization categories are denoted by $\mathcal{C}^W$. \par The local cohomology of $Y \subset X$ is $H_Y (X)$ and the Schur functor is $\BBS^{\alpha}$. \par The following are specific key sets, spaces or categories, defined here as a reference point, as well as where they are first defined. 
\begin{itemize}
	\item $\mathcal{UP}_{n,k}$ - Strictly upper triangular Young diagrams - p16
	\item $\Lambda = \End(\medoplus_{\alpha \in \mathcal{UP}_{n,k}} (\BBS^{\alpha}S^{*}\medotimes \Sym{\bullet}X^*)^G  )$ - p22
	\item $\mathcal{X} = \left[\Hom(S,V)/SL(S)\right]$ - p38
	\item $\mathcal{Z} = \left[\Hom(S,V) \times \wedge^k S /GL(S) \right]$ - p38
	\item $\mathcal{X}_L = \left[\Hom(S,V)\times (\wedge^k S\setminus 0) \times (L\setminus 0)/GL(S)\times \BBC^*_{0,-1,1} \right]$ - p44
	\item $\mathcal{Z}_L = \left[\Hom(S,V)\times \wedge^k S \times (L\setminus 0)/GL(S)\times \BBC^*_{0,-1,1} \right]$ - p44
	\item $Y_L = \left[\Hom(S,V)^f\times \wedge^k S \times (L\setminus 0)/GL(S)\times \BBC^*_{0,-1,1} \right]$ - p44
	\item $\mathcal{B}^W_L \subset D^b(Y_L, W)$ - p43
	\item $\mathcal{D}^W_L \subset D^b(\mathcal{X}_L, W)$ - p45
\end{itemize}

\section{Background}\label{sec:background}
In this section we will briefly cover some of the technical results we use as well as providing a very brief introduction to some of the possibly less well known theories that we use.
\subsection{Schur functors and Young diagrams}\label{sec:schurandyoung}
We mainly work with vector bundles associated to representations, Young diagrams and Schur functors are a useful way to work with such objects.
\begin{definition}
	A {\it Young diagram}, $\alpha = (\alpha_{1}, \alpha_{2}, \dots , \alpha_{n})$, is an element of $\BBZ^n_{\geq 0}$ subject to the condition $\alpha_{1} \geq \alpha_{2} \geq \dots \geq \alpha_{n} \geq 0$.\\
	We call this a {\it Young diagram of size $n$}. If in addition, $\alpha_{1} = k$, we call it a {\it Young diagram of size $n \times k$}.
\end{definition}
\noindent
We can represent $\alpha$ as a diagram of boxes with $n$ rows, the $i^{th}$ row has $\alpha_{i}$ boxes, see Figure 1 for an example.
\begin{figure}[h!]
	\centering
	\includegraphics[width=.6\linewidth]{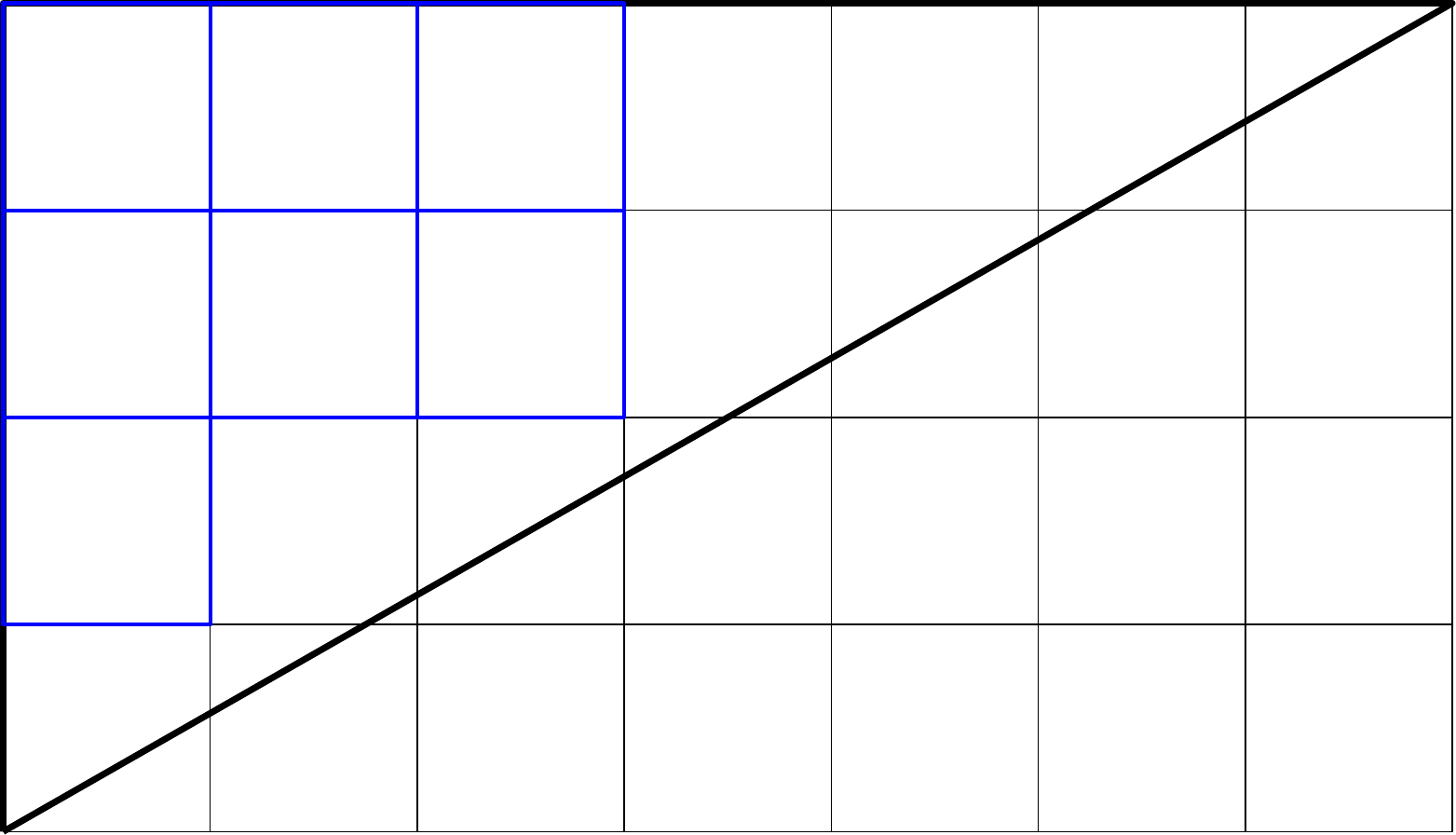}
	\caption{The Young diagram  $(3,3,1,0)$}
\end{figure}
\newline
One can classify $SL$ and $GL$ representations using Young diagrams. We have
$$\{\mbox{Irreducible}\; SL_n \; \mbox{representations}\} \cong \{\mbox{Young diagrams of size}\; n-1\}.$$
Often we will use Young diagrams of size $n$ with $\alpha_{n}=0$ for describing $SL_n$ representations. To describe $GL$ representations we will extend Young diagrams by dropping the positivity condition.
$$\{\mbox{Irreducible}\; GL_n \; \mbox{representations}\} \cong \{\mbox{Extended Young diagrams of size}\; n\}.$$
To go from $\alpha$, a Young diagram, to a representation we use Schur functors.
\begin{definition}
	The \emph{Schur functor}, $\BBS$, takes Young diagrams of size $n$ to vector spaces.\\
	Let $\alpha$ be a Young diagram of size $n$, we write $\BBS^{\alpha} \BBC^n$ for the associated vector space.
\end{definition}
\noindent
The actual definition of $\BBS^{\alpha} \BBC^n$ is technical, we will not provide it here, see \cite{fulton2013representation}. All we need to know is:
\begin{itemize}
	\item $\BBS^{(i, 0 \dots , 0)} V \cong \Sym{i} V$.
	\item $\BBS^{(i, i \dots , i)} V \cong \left(\det V\right)^i$.
	\item $\BBS^{\alpha} V$ is the irreducible representation of $SL(V)$ with highest weight $\alpha$.
\end{itemize}
We use the final statement for extending the definition of $\BBS$ to extended Young diagrams. More accurately, if $\alpha_{n} < 0$, then we can define $\BBS^{\alpha} V = \BBS^{(\alpha_{1} - \alpha_{n}, \alpha_{2} - \alpha_{n}, \dots , \alpha_{n} - \alpha_n)}V \medotimes \left(\det V^*\right)^{-\alpha_{n}}$.\par
One of the most useful things about describing representations this way is we can calculate what representations appear in $\BBS^{\alpha} V \medotimes \BBS^{\beta} V$ using the Littlewood-Richardson rule for Young diagrams. Again we will not provide the exact details of this rule, see \cite{fulton2013representation} for more information. Note however that if $\BBS^{\gamma}V \subset \BBS^\alpha V\medotimes \BBS^{\beta}V$ then the number of boxes in $\gamma$ is the sum of the boxes in $\alpha$ and $\beta$. We will explain additional details of the rule when we require them.\par
A very important set of Young diagrams for us is 
$$\mathcal{UP}_{n,k}  =\{\mbox{Strictly upper triangular Young diagrams of size at most}\; k \times (n-k) \}.$$
Algebraically we have $\alpha \in \mathcal{UP}_{n,k}$ if and only if $\alpha_i \leq (k-i)(n-k)/k$. These give us a specific collection of $SL_k$ or  $GL_k$ representations. It is this collection that will be used for the NCCR and this collection also defines the first block of the minimal Lefschetz decomposition of $Gr(k,n)$. See Figure 1 for an example of such a Young diagram in $\mathcal{UP}_{11,4}$.
\subsection{Borel-Weil-Bott}\label{sec:borel_weil_bott}
The Borel-Weil-Bott theorem is used to calculate cohomology of vector bundles that correspond to representations. Almost all the sheaves we work with fit into this class.
\begin{theorem}[Borel-Weil-Bott]
	Let $G$ be a reductive group, let $B \subset G$ be a Borel subgroup. Let $\alpha$ be a weight of $B$, it corresponds to a line bundle $\mathcal{O}(\alpha)$ on $G/B$. Let $\rho$ be half the sum of the positive roots of $G$. There exists a unique element, $\sigma$, of the Weyl group of $G$ that takes $\alpha + \rho$ to a dominant weight, $\sigma(\alpha + \rho)$. Let $\alpha^+ = \sigma(\alpha + \rho) - \rho$. Then we have two possible situations.
	\begin{itemize}
		\item $H^i\left( G/B, \mathcal{O}(\alpha)\right) = 0$ for all $i$ when $\alpha^+$ is not dominant.
		\item $H^i \left( G/B, \mathcal{O}(\alpha)\right) = V(\alpha^+)^*$ for $i=l(\alpha^+)$ and zero for all other $i$ when $\alpha^+$ is dominant.
	\end{itemize}
	Here $V(\alpha^+)^*$ denotes the dual of the irreducible representation of $G$ with highest weight $\alpha^+$ and $l(\alpha^+)=l(\sigma)$ is the smallest number $k$ such that $\sigma$ is the product of $k$ simple reflections, called the {\it length} of $\sigma$.
\end{theorem}
\begin{remark}
	$G/B$ is a flag variety, in the special case $G=GL(V)$, we have $G/B$ is the complete flag on $V$.
\end{remark}
\noindent
We will frequently need to calculate cohomology on the Grassmannian. We can do this using Borel-Weil-Bott. \par
Let $V$ be a vector space of dimension $n$. Consider $Gr(k,V) = Gr(V,n-k)$, the Grassmannian of $k$-planes in $V$, or identically the same as $(n-k)$-quotients of $V$. We have the tautological short exact sequence
$$0 \rightarrow S \rightarrow V \medotimes \mathcal{O}_{Gr(k,V)} \rightarrow Q \rightarrow 0 $$
where $S$ is the tautological sub-bundle and $Q$ is the tautological quotient bundle. 

\begin{proposition}\label{prop:borel_weil_bott}
	Let $\beta \in \BBZ^k, \gamma \in \BBZ^{n-k}$ be two non-increasing sequences of integers and let $\alpha = (\beta, \gamma) \in \BBZ^n$. Then we have
	$$ H^{\bullet}(Gr(k,V), \BBS^{\beta} S^* \medotimes \BBS^{\gamma} Q^*) \cong \BBS^{\alpha'} V^*[-l(\alpha')]$$
	where $\alpha'$ is the unique dominant weight in the twisted Weyl group orbit of $\alpha$, and $l(\alpha')$ is the length. (If such an element does not exist, the cohomology is zero.) 
\end{proposition}
This proposition is the result we will mainly use. It follows from applying Borel-Weil-Bott twice. First we apply it to the two flag bundles, $F(S)$ and $F(Q)$ over $Gr(k,V)$, then by using the K{\"u}nneth formula we find out that $\mathcal{O}(\alpha)$ on $F(S)\times F(Q)=F(V)$ pushes down to $\BBS^{\beta} S^* \medotimes \BBS^{\gamma} Q^*$. Second we apply it to $G=GL(V)$ with $G/B = F(V)$ and the line bundle $\mathcal{O}(\alpha)$.\par
As a special case, if $\alpha$ is dominant then we have $H^0$ is $\BBS^{\alpha} V^*$ and there is no higher cohomology. Also note, by considering the vector space $V^*$ instead, any result we get for $S^*, Q^*$ also holds for $S,Q$ up to appropriate swapping of $k$ and $n-k$.
\subsection{Local Cohomology}\label{sec:localcoh}
Let $X$ be a scheme, $\mathcal{F}$ a sheaf on $X$ and let $Y \subset X$ be a closed subscheme, we can then define $\Gamma_Y(\mathcal{F})$ to be the sections of $\mathcal{F}$ with support in $Y$, deriving we get $H^i_Y(X,\mathcal{F})$, which are called the {\it local cohomology groups} of $Y$ in $X$. There is also another definition, let $U = X \setminus Y$ and let $j: U \rightarrow X$ be the inclusion, then $\mathcal{H}^{\bullet}_Y(X,\mathcal{F}) = Cone(\mathcal{F} \rightarrow j_* j^* \mathcal{F})$ where the map comes from adjunction. Take global sections to get $H^i_Y(X,\mathcal{F})$. Note that the projection formula gives us $\mathcal{H}^{\bullet}_Y(X,\mathcal{F}) = \mathcal{H}^{\bullet}_Y(X,\mathcal{O}_X)\medotimes \mathcal{F}$.\par
There are two sequences that hold for local cohomology that will be useful to us. First we have the following long exact sequence.
\begin{align*}
0 \rightarrow &H^0_Y(X, \mathcal{F}) \rightarrow H^0(X,\mathcal{F}) \rightarrow H^0(U, \mathcal{F}|_U) \rightarrow \\
&H^1_Y(X, \mathcal{F}) \rightarrow H^1(X, \mathcal{F}) \rightarrow H^1(U, \mathcal{F}|_U) \rightarrow \cdots
\end{align*}
For the second sequence, let $Z \subset Y \subset X$ be a sequence of closed subspaces, we then have an exact triangle
\begin{equation*}
	\mathcal{H}^{\bullet}_Z(X, \mathcal{F}) \rightarrow \mathcal{H}^{\bullet}_Y(X, \mathcal{F}) \rightarrow \mathcal{H}^{\bullet}_{Y\setminus Z}(X \setminus Z, \mathcal{F}|_{X\setminus Z}).
\end{equation*}
We will often use local cohomology in the following way. Let $X$ be a scheme and $\mathcal{M}$ be a bundle on $X$, then given a family $Z$ of closed subschemes over $X$ that lives inside $\mathcal{M}$ we have $\mathcal{H}^i_Z(\mathcal{M}, -)$, we will use this to denote the pushdown to $X$. In the special case of $\mathcal{N} \subset \mathcal{M}$, where $\mathcal{N}$ is a sub-bundle of $\mathcal{M}$ we know exactly what the local cohomology of $\mathcal{O}_{\mathcal{M}}$ is.
\begin{lemma}[{\cite[Prop. 3.3.1]{van1993cohen}}]\label{lem:localcoh_for_subbundles}
	Let $X$ be a smooth scheme, let $\mathcal{M}$ be a bundle over $X$ and $\mathcal{N}$ be a sub-bundle of $\mathcal{M}$. Then
	$$\mathcal{H}^{\codim \; \mathcal{N}}_{\mathcal{N}} (\mathcal{M}, \mathcal{O}_{\mathcal{M}}) = \Sym{\bullet}\left(\mathcal{N}^*\right) \medotimes \Sym{\bullet} \left(\mathcal{M}/ \mathcal{N}\right)\medotimes \det \left(\mathcal{M}/ \mathcal{N}\right)$$
	and all other cohomology is zero.
\end{lemma}
\noindent
Finally there is a third way of thinking about local cohomology, it is the same as \emph{relative De Rham homology}, which can be defined as follows. Given a diagram 
\begin{center}
	\begin{tikzcd}
		Y \arrow[rd] \arrow[r, hook] & \tilde{X} \arrow[d] \\
		& X
	\end{tikzcd}
\end{center}
where the horizontal map is a closed immersion and the vertical one is smooth, let $H^{DR}_i(Y/X)$ be defined as $$H^{DR}_i(Y/X) = H^{2n-i}_Y(\tilde{X},\Omega^{\bullet}_{\tilde{X}/X})$$ where $n$ is the relative dimension of $\tilde{X}/X$. It is functorial for proper maps, so if we have $Y' \rightarrow Y$ proper, we get a map $H^{DR}_i(Y'/X) \rightarrow H^{DR}_i(Y/X)$. We also have 
$$ H^{DR}_{-i}(Y/X) = H^i_Y(X, \mathcal{O}_X)$$
for $Y \subset X$ closed. For more details on relative De Rham homology see \cite{van1991cohen}.
\subsection{Semi-orthogonal Decompositions and Admissible Subcategories}\label{sec:decompandadmissible}
Let $\mathcal{A}, \mathcal{B} \subset \mathcal{C}$ be triangulated subcategories of a triangulated category. If we have $\Hom (\mathcal{B}, \mathcal{A}) = 0$, and for all $X \in \mathcal{C}$, we can find an exact triangle $B \rightarrow X \rightarrow A$ where $A \in \mathcal{A}, B \in \mathcal{B}$, we say that $\mathcal{C}$ has a {\it semi-orthogonal decomposition} and we write $\mathcal{C} = \langle \mathcal{A}, \mathcal{B} \rangle$. This extends naturally to semi-orthogonal decompositions of longer, or even infinite length.\par
A source of semi-orthogonal decompositions are admissible subcategories. 
\begin{definition}
	Let $\mathcal{A}$ and $\mathcal{C}$ be triangulated subcategories, and let $i: \mathcal{A} \rightarrow \mathcal{C}$ be fully faithful.
	\begin{itemize}
		\item We call $\mathcal{A}$ {\it left (right) admissible} if $i$ has a left (right) adjoint. 
		\item We call $\mathcal{A}$ {\it admissible} if it is left and right admissible.
	\end{itemize}
\end{definition}
\begin{lemma}
	Let $\mathcal{A} \subset \mathcal{C}$ be left admissible, then $\mathcal{C} = \langle \mathcal{A}, {}^{\perp}\mathcal{A} \rangle$. Let $\mathcal{A} \subset \mathcal{C}$ be right admissible, then $\mathcal{C} = \langle \mathcal{A}^{\perp}, \mathcal{A} \rangle$.
\end{lemma}
\noindent
Here, ${}^\perp \mathcal{A}$ is the \emph{left orthogonal} to $\mathcal{A}$, it is all objects $X \in \mathcal{C}$ such that $\Hom(X, \mathcal{A}) = 0$, and similarly, $\mathcal{A}^{\perp}$ is the \emph{right orthogonal}.\par
If we have a bounded derived category of a scheme or algebra then one way to show a subcategory is admissible is as follows. First find a finite generating set and then show that the endomorphism algebra of the generators has finite global dimension, as then $\Hom(T, -)$ and $T \medotimes_{\End T} -$ where $T$ is the direct sum of the generators, give the needed functors. In particular, this shows that any subcategory that is generated by exceptional objects is admissible.
\subsection{Global Quotient Stacks and Matrix Factorizations}\label{sec:stackandmf}
Let $X$ be a scheme and $G$ be a reductive group acting on $X$, we can then form the global quotient stack $\left[X/G\right]$. These kinds of objects will appear often, and in most cases $X$ will be a vector space, or an open set in a vector space. For practical purposes we can treat these as standard schemes. Two simple results to help with understanding are:
\begin{enumerate}
	\item $\mbox{coh} \left[X/G\right] = \mbox{coh}_G X.$
	\item If $G$ acts freely on $X$ then the stack is isomorphic to the quotient scheme, i.e. let $\BBC^*$ act with weight 1 on $\BBC^n$, then $\left[ (\BBC^n \setminus 0)/ \BBC^*\right] \cong \BBP^{n-1}$.
\end{enumerate}
We will frequently use a generalization of $D^b(X)$, known as matrix factorizations. The idea for this comes from algebra, functions that can not be factored in the standard way can often be factored by matrices. Let $R$ be a commutative ring and let $f \in R$. Then a matrix factorization of $f$ is a pair of free modules $P,Q$ and maps $\alpha: P \rightarrow Q, \beta: Q \rightarrow P$ such that $\alpha \circ \beta = f Id_Q$ and $\beta \circ \alpha = f Id_P$. We want to generalize this and put it into a triangulated framework. \par
Let $X$ be a scheme/stack and let $W:X \rightarrow \BBC$ be a {\it superpotential} (a global section of the structure sheaf). Then we could ask for a pair of sheaves and maps between them such that both compositions equal acting by $W$ on the sheaves. Instead we will do something slightly different.\par
We add a $\BBC^*$-action to $X$, called the $R$-charge, this gives us a shift operation, $[1]$, given by tensoring with a character of $\BBC^*$. We pick the $R$-charge such that $W$ has weight 2. Then we can consider sheaves $\mathcal{F}$, which have a map $d: \mathcal{F} \rightarrow \mathcal{F}[1]$ such that $d^2 = W Id_{\mathcal{F}}$. If we set this up formally we end up with a triangulated category $D^b(X, W)$ with objects $(\mathcal{F}, d)$, called matrix factorizations. We choose not to denote the $R-$charge in this notation. 
\begin{remark}
	If $X$ itself is of the form $X'/G$ then we can modify the $R-$charge by any central one parameter subgroup of $G$ without changing $D^b(X, W)$, see \cite[Remark 2.5]{2017arXiv170501437V}.
\end{remark}
\begin{remark}
	If the structure sheaf exists only in even degrees then we can split $\mathcal{F}$ into odd and even parts and we end up with a situation looking very much like the original motivation.
\end{remark}
\begin{remark}
	Note that this is a generalization of the original derived category, as if $W=0$, then setting the $R-$charge to be trivial, we recover the definition of $D^b(X)$, i.e. $D^b(X,0) \cong D^b(X)$, see \cite[Prop. 2.1.6]{2013arXiv1306.3957B}	
\end{remark}
An important result is called Kn{\" o}rrer periodicity and it relates derived categories of subschemes with matrix factorization categories of line bundles.\par
Let $X$ be a scheme, and $\mathcal{L}$ be a line bundle on $X$, then given a section $f$ of $\mathcal{L}^*$ we can define a superpotential on the total space of $\mathcal{L}$ by $pf$ where $p$ is the fibre coordinate on $\mathcal{L}$. We can also define an $R-$charge by letting $\BBC^*$ act with weight 2 on the fibres and trivially on $X$. Let $Y \subset X$ be the zero locus of $f$. Then we have the following result.
\begin{theorem}[Kn{\" o}rrer periodicity, \cite{Shipman_2012}, \cite{Hirano_2017} ]	
	$D^b(Tot(\mathcal{L}), pf) \cong D^b(Y)$.
\end{theorem}
This result in fact holds more generally. Let $\mathcal{F}$ be a vector bundle and $f$ a section of the dual vector bundle, then $D^b(Tot(\mathcal{F}), pf) \cong D^b(Z(f))$ as long as $Z(f)$ has the right dimension. See \cite{Shipman_2012}, \cite{Hirano_2017} for the exact technical details. \par
We also commonly use a specific type of subcategory of $D^b(X,W)$. Let $A$ be a set of vector bundles on a stack $\mathcal{X}$ and let $W$ be a superpotential on $\mathcal{X}$. Then we define the subcategory $D^b(A, W) \subset D^b(\mathcal{X}, W)$ as all matrix factorizations quasi-isomorphic to matrix factorizations built from vector bundles in $A$. The reason for this definition is we want to study matrix factorizations built from a finite set of vector bundles, but these are not closed under isomorphism.\par
An useful result to help study these subcategories is:
\begin{lemma}\label{lem:equivalentdescriptionofwindowsubcategory}
	Let $\mathcal{V} = [V/G]$ for $V$ a vector space, such that there is a one-parameter subgroup of $G \times \BBC^*$ which makes the structure sheaf non-negatively graded. (The $\BBC^*$ is the $R$-charge.) Let $\mathcal{A}$ be a set of irreducible representations of $G$ and let $A$ be the set of associated vector bundles. Then $\mathcal{E} \in D^b(A,W) \subset D^b(\mathcal{V}, W)$ if and only if the homology of $\mathcal{E}|_0$ contains only irreps from $\mathcal{A}$.
\end{lemma}
\begin{proof}
	See \cite[Lemma 2.8]{2016arXiv160904045V} and the remark after that Lemma.
\end{proof}
Informally this result says that the condition of asking a matrix factorization to be built from a set of vector bundles associated to representations is equivalent to asking that matrix factorization restricted to the origin only has those representations appearing. The second condition is more general and mathematically a better definition, the first condition can be easier to work with. If we are in this situation, we can then equivalently define $D^b(A, W)$ as $\{\mathcal{E} | \mathcal{E}|_0 \in \mathcal{A} \} \subset D^b(\mathcal{X}, W)$.\par
For technical details on global quotient stacks and matrix factorizations see \cite{2012arXiv1203.6643B}, \cite{2014arXiv1401.3661A}.

\section{Algebraic NCCR}\label{sec:algebraicnccr}
The first thing we will do is find an NCCR for the affine cone of the Grassmannian.\par
Let $S$ and $V$ be vector spaces of dimension $k$ and $n$. Assume $1 < k  < n - 1$ and that $k$ and $n$ are coprime\footnote{If $k=1$ or $n-1$ then there is no singularity.}. Also let $X=\Hom(S, V)$, then $G=SL(S)$ acts on $X$. We want to find a NCCR for $X/G$. There are no general methods for finding NC(C)R's of an arbitrary singular scheme, however for specific types of singularities there are some approaches. Our singularity is a quotient by a reductive group, for this type of singularity, {\v S}penko and Van den Bergh \cite{2015arXiv150205240S} provide a general method for finding NC(C)R's.\par
Based on the finite group case, for $H$ a reductive group acting on $Y$, they consider the algebra
$$\Lambda_{\mathcal{L}} = \left(\End(U)\medotimes k[Y]\right)^H$$
where $U = \medoplus_{\alpha \in \mathcal{L}} \BBS^{\alpha}$ and $\BBS^{\alpha}$ is the irreducible representation of $H$ with highest weight $\alpha$, $\mathcal{L}$ is any finite collection of highest weights. Note that in general the algebra has a natural grading as it can be represented by a quiver algebra and in our case, the singularity has a natural $\BBC^*$ action given by dilation.\par
We want this graded algebra to have finite global dimension and to be Cohen-Macaulay. It splits into components $\left(\BBS^{\beta} \medotimes k[Y]\right)^H$ where $\BBS^{\beta}$ appears in the decomposition of $(\BBS^{\alpha})^* \medotimes \BBS^{\alpha'}$ for $\alpha, \alpha' \in \mathcal{L}$ and it is sufficient to prove that each individual component is Cohen-Macaulay. {\v S}penko and Van den Bergh reference a result from \cite{van1991cohen} to show Cohen-Macaulayness and we will strengthen this result in our situation, the majority of this section is dedicated to doing this.\par
Most of the paper \cite{2015arXiv150205240S} is dedicated to showing finite global dimension, they first show that is sufficient to prove that each $H$-equivariant $\Lambda$-module of the form \\
$P_{\alpha} = \Hom \left(U\medotimes k[Y], \BBS^{\alpha} \medotimes k[Y]\right)$ has finite projective dimension. Then they build a resolution of $P_{\alpha}$ as follows: Pick a one parameter subgroup $\lambda$ of $T \subset B \subset H$ ($T$ a maximal torus inside a Borel subgroup, $B$). Let $Z_{\lambda} \subset Y$ be the subspace of points that flow to zero as $\lambda$ tends to zero. Let $\widetilde{Z_{\lambda}} = Z_{\lambda}\times_B H$, it is a resolution of $H Z_{\lambda}$ and a sub-bundle of $\widetilde{Y} = Y \times_B H$ where both are bundles over $H/B$. Taking the Koszul resolution of $\widetilde{Z_{\lambda}}$ gives an exact complex, we can then tensor it with $\BBS^{\alpha}$ and then applying $\Hom(U\medotimes k[Y],-)$ gives us a complex of $\Lambda$-modules. If $\alpha$ and $\lambda$ are suitably related then the cohomology of this complex is exact and is a resolution of $P_{\alpha}$. They also describe all the modules $P_{\beta}$ that appear in this resolution.\par
The idea to show finite global dimension is then as follows: Take any module $P_{\alpha}$, if $\alpha \in \mathcal{L}$, then $P_{\alpha}$ is projective. If $\alpha \not \in \mathcal{L}$ then we use the resolution of $P_{\alpha}$ constructed above, the aim is to prove that by iterating these resolutions a finite number of times we can find a projective resolution of $P_{\alpha}$. If we can do this for all $\alpha$, we have finite global dimension, as the $P_{\alpha}$ are all you need to resolve the graded simples. We will explain briefly in Section \ref{sec:fgdim} how Fonarev, \cite{fonarev2013minimal} shows this. \par
Note that an Auslander-Buchsbaum type equation, \cite[Lemma 2.16]{2010arXiv1007.1296I}, tells us that once we have a NCCR, $\Lambda$, the global dimension of $\Lambda$ will be equal to $\dim Y/H$. So even though we only ask for finite global dimension, in fact, once we have Cohen-Macaulayness, $\Lambda$ has the right dimension to be considered a resolution. 
\subsection{Setup}\label{sec:setup}
Recall, $\mathcal{UP}_{n,k}$ is the set of all Young diagrams $\alpha$ such that $\alpha_{i} \leq (n-k)\frac{k-i}{k}$, $X = \Hom(S,V), G=SL(S)$, we want to show that the graded algebra 
$$\Lambda = \End \left(\bigoplus_{\alpha \in \mathcal{UP}_{n,k}} (\BBS^{\alpha}S^{*}\medotimes \Sym{\bullet}X^*)^G \right)$$
gives a NCCR of $X/G$. Note that another way of describing this result is the fact that the subcategory
$$\mathcal{D} = \langle \mathcal{O}_{X} \medotimes \BBS^{\alpha}S^* \rangle_{\alpha \in \mathcal{UP}_{n,k}} \subset D^b\left([\Hom(S, V)/SL(S)]\right)$$
is a crepant categorical resolution as $\mathcal{D} \simeq D^b(\Lambda)$ once we have finite global dimension. \par
We will do this by using the framework of {\v S}penko and Van den Bergh. First we show that $\Lambda$ is Cohen-Macaulay, to do that we need to show that each representation that appears in $\Lambda$ is Cohen-Macaulay. We use the techniques from \cite{van1991cohen} to find all Cohen-Macaulay modules of covariants. The same result is proved using different methods in \cite{raicu2014local}. Second we will show that $\Lambda$ has finite global dimension in Section \ref{sec:fgdim}.\par 
To show Cohen-Macaulayness of $\Lambda$ we need to show Cohen-Macaulayness for each component, which are of the form $\left(\BBS^{\alpha}S^{*}\medotimes \Sym{\bullet} X^*\right)^G$. Now let $R_{\alpha}^G$ be the component of $\mathcal{O}_X(X)$ with weight $\alpha$ (split $\mathcal{O}_X(X)$ into irreducible representations of $G$). We then have $R_{\alpha}^G$ is Cohen-Macaulay if and only if $\left(\BBS^{\alpha}S^{*}\medotimes \Sym{\bullet} X^*\right)^G$ is Cohen-Macaulay, \cite[Lemma 3.2]{van1989trace}. Therefore it is sufficient to show $R_{\alpha}^G$ is Cohen-Macaulay and we can do this by using the following result.
\begin{lemma}[{\cite[Cor. 4.2]{van1989trace}}]\label{lem:localcoh_and_cm}
	Let $X^u \subset X$ be the null cone, let $h = \dim \mathcal{O}_X^G$, then $R_{\alpha}^G$ is Cohen-Macaulay if and only if there is no representation with weight $\alpha$ in $H_{X^u}^i ( X, \mathcal{O}_X)$, where $i=0,1,\dots ,h-1$.
\end{lemma}
\noindent
As $X^u$ is singular, it is hard to calculate $H_{X^u}^i ( X, \mathcal{O}_X)$ directly so we will use a spectral sequence from \cite{van1991cohen} that converges to $H_{X^u}^i ( X, \mathcal{O}_X)$. \par
The actual calculations turn out to be quite technical, so we split the work into two parts, first in Section \ref{sec:SS}, we give a geometric motivation for the spectral sequence and describe the terms. Secondly, in Section \ref{sec:calculations}, we do all the calculations that are needed. The result we get is Lemma \ref{lem:representationsappearinginlocalcohomology} and then proving $\Lambda$ is Cohen-Macaulay, Proposition \ref{prop:algebraisCM}, follows directly. The rest of Section \ref{sec:calculations} can be passed over without losing any understanding as it is just cohomology calculations on the varieties described in Section \ref{sec:SS}.\par
Fix a Borel subgroup $B \subset G$, all parabolic subgroups will contain $B$. Let $w_1 = (1,0, \dots, 0), w_2 = (1,1,0 \dots , 0), \dots, w_{k-1} = (1,\dots ,1,0)$.
These are the fundamental weights of $SL(S)$, i.e. any dominant weight can be written as $\sum a_i w_i $ where $a_i \geq 0$. Note, the weights of $SL(S)$ live in $\BBR^k$ up to shifts of $(a,a, \dots , a)$, so we shift all our weights to have a zero in the final position.
\subsection{The Spectral Sequence}\label{sec:SS}
We want to calculate $H_{X^u}^i ( X, \mathcal{O}_X)$, to do this we will use a spectral sequence where the terms are the local cohomology groups $H^*_Y(\tilde{X}, -)$ for $Y$ a sub-bundle of $\tilde{X} = X \times G/P$, for $P$ some parabolic subgroup. $Y$ will be closely related to a resolution of a stratification of $X^u$. \par
We will do the first two cases by hand, then give the general construction in \cite{van1991cohen} applied in our situation. Let $\dim V = 2$, then $X^u = \{rank \leq 1 \;\, maps\}$ and $0 \subset X^u \subset X$, so we get an exact triangle
$$H^{\bullet}_0(X, \mathcal{O}_X) \rightarrow H^{\bullet}_{X^u}(X, \mathcal{O}_X) \rightarrow H^{\bullet}_{X^u \setminus 0}(X\setminus 0 , \mathcal{O}_{X\setminus 0}).$$
We think of this as replacing $H^{\bullet}_{X^u}(X, \mathcal{O}_X)$ with $H^{\bullet}_0(X, \mathcal{O}_X)$, which we understand, and \\$H^{\bullet}_{X^u \setminus 0}(X\setminus 0 , \mathcal{O}_{X\setminus 0})$. To deal with this term we use the relation between local cohomology and relative De Rham homology. We have the following diagram
\begin{center}
	\begin{tikzcd}
		X^u\setminus 0 \arrow[d, hook] \arrow[rrd] \arrow[r, "\sim"] & {Y\setminus Gr(S,L)} \arrow[r, hook] & {(X \setminus 0) \times Gr(S,L)} \arrow[d] \\
		X \setminus 0 \arrow[rr] &  & X \setminus 0
	\end{tikzcd}
\end{center}
where $Y = \Hom(L, V)_{Gr(S,L)}$, $\dim L = 1$. $Y$ is a resolution of $X^u$ and also a sub-bundle of the trivial bundle over $Gr(S,L)$ with fibres $X \setminus 0$. This diagram shows that we have
$$H^{\bullet}_{X^u\setminus 0}(X\setminus 0, \mathcal{O}_{X\setminus 0}) \cong H^{DR}_{\bullet}((X^u\setminus 0 )/(X \setminus 0)) \cong H^{\bullet}_{Y\setminus Gr(S,L)}((X \setminus 0) \times Gr(S,L), \Omega^{\bullet}_{(X \setminus 0) \times Gr(S,L)/(X \setminus 0)}).$$
Then using the exact triangle for local cohomology, we can replace \\ $H^{\bullet}_{Y\setminus Gr(S,L)}((X \setminus 0) \times Gr(S,L), \Omega^{\bullet}_{(X \setminus 0) \times Gr(S,L)/(X \setminus 0)})$ with $H^{\bullet}_{Y}(X \times Gr(S,L), \Omega^{\bullet}_{X \times Gr(S,L)/X})$ and \\
$H^{\bullet}_{0 \times Gr(S,L)}(X \times Gr(S,L), \Omega^{\bullet}_{X \times Gr(S,L)/X})$. We can calculate both of these using Lemma \ref{lem:localcoh_for_subbundles}, so we are done. (Up to working out maps and gradings.) \par
The next case is $\dim V = 3$, here we have $0 \subset X^1 \subset X^2 = X^u$, where $X^i$ are the maps of rank $\leq i$. If we consider $0 \subset X^2 \subset X$ as above we get the following. Let $Y=\Hom (L, V)_{Gr(S,L)}$ where $\dim L = 2$. Again $Y$ is a resolution of $X^u$, however this time $X^u\setminus 0 \not \cong Y \setminus Gr(S,L)$. Over the rank two maps the map $Y\rightarrow X^u$ is an isomorphism, however over the rank one maps we have a fibre isomorphic to $\BBP^1$. What is instead true is that $X^u \setminus X^1 \cong Y\setminus Y^1$ where $Y^1$ is the preimage of $X^1$. \par
This means that we should consider $0 \subset X^1 \subset X$ and $X^1 \subset X^2 \subset X$ and in the same way as above, we can replace the local cohomology of $X^u$ by local cohomology of $0$, $Y$, $0\times Gr(S,L')$, $Z= \Hom(L', V)_{Gr(S,L')}$ (where $\dim L'$ = 1) and a term $Y^1$. We can also deal with this final term in a similar way. We have $U= \Hom(M, V)_{Fl(3,2,1)}$ resolving $Y^1$ ($M$ is rank 1) and $Y^1\setminus 0\times Gr(S,L) \cong U \setminus 0 \times Fl(3,2,1)$ so we can replace the local cohomology of $Y^1$ by that of $0\times Gr(S,L)$, $U$ and $0\times Fl(3,2,1)$. All of these terms are linear sub-bundles of trivial bundles so we can use Lemma \ref{lem:localcoh_for_subbundles} to calculate their local cohomology.\par
They all fit together as in the following diagram
\begin{center}
	\begin{tikzcd}
		& {0 \times Gr(S,L)} \arrow[ld] \arrow[dd] &                                                     & {0 \times Fl(3,2,1)} \arrow[ld] \arrow[ll] \arrow[dd] \\
		{\Hom (L, V)_{Gr(S,L)}} \arrow[dd] &                                          & {\Hom(M, V)_{Fl(3,2,1)}} \arrow[ll] \arrow[dd] &                                                       \\
		& 0  \arrow[ld]                            &                                                     & {0 \times Gr(S,L')} \arrow[ll] \arrow[ld]              \\
		X^u                                     &                                          & {\Hom(L', V)_{Gr(S,L')}} \arrow[ll]              &                                                      
	\end{tikzcd}
\end{center}
In general we have exactly the same idea. Stratify $X^u$ by $X^i = \{rank \leq i \;\, maps\}$, take a specific resolution, $Y^i$, of each $X^i$ and then consider a stratification of each $Y^i$ given by the preimages of the $X^j$ with $j<i$. Resolve those stratifications and repeat this procedure. The general case is explained in \cite{van1991cohen}, here is a very brief summary of the result in our situation using some of their language.\par
Consider the stratification $0 \subset X^1 \subset \cdots \subset X^{k-1} = X^u$ where $\dim S = k$.\\
Let $Y^i$ be a linear subspace of $X$ such that $G Y^i = X^i$. Let $P \supseteq B$ be any parabolic subset preserving $Y^i$.  We get a spectral sequence whose terms are relative De Rham homology.
$$ \medoplus_{(i,P) \in \mathcal{R}} H^{DR}_{\bullet}(Y^i \times_P G / X) \implies H^{\bullet}_{X^u}(X, \mathcal{O}_X)$$
For the grading, definition of $\mathcal{R}$ and the maps see \cite[Theorem 5.2.1]{van1991cohen}. We can replace $H^{DR}_{\bullet}(Y^i \times_P G / X)$ with a simple local cohomology calculation as $Y^i \times_P G$ is a subbundle of  $ X \times_P G \cong X \times (G/P)$, this is as $X$ has a $G$ action and so we get an isomorphism $(x, g) \mapsto (xg, [g])$. Therefore letting $Y$ be any of the terms $Y^i \times _P G$, once we have calculated $H^{\bullet}_{Y}(X \times G/P)$, we will know what terms can not appear in $H^{\bullet}_{X^u}(X, \mathcal{O}_X)$.
\subsection{Calculations}\label{sec:calculations}
Now that we know what terms appear in the spectral sequence we next calculate all the possible representations that can appear in these terms, this will also tell us what representations can not appear in $H^{\bullet}_{X^u}(X, \mathcal{O}_X)$. The final result is
\begin{lemma}\label{lem:representationsappearinginlocalcohomology}
	Any representation appearing in $H_{X^u}^i ( X, \mathcal{O}_X)$ for $i<h$ has weight $\sum a_i w_i$ where at least one $a_i\geq n-k$ and all the $a_i$ are non-negative.
\end{lemma}
This is strong enough to prove that $\Lambda$ is a Cohen-Macaulay algebra, see Proposition \ref{prop:algebraisCM} for a proof. The rest of this section will be devoted to proving Lemma \ref{lem:representationsappearinginlocalcohomology}, it is mainly just calculations and can be skipped, the geometric picture appeared in the previous section. A slightly stronger result (an if and only if statement) was proved independently using another method of calculating local cohomology, \cite[Theorem 4.6]{raicu2014local}.\par
Lemma \ref{lem:localcoh_and_cm} tells us that we only need to check what representations appear in the first $h = \dim X - \dim G$ degrees. We have $\dim X = \dim \Hom(S, V) = kn$ and $\dim G = \dim SL(S) = k^2-1$ which gives us that $h = nk - k^2 + 1$. \par
Let $F = G/P$ for any parabolic subgroup $P$, it is a partial flag variety for $V$. We have $H^{\bullet}_Y(X', \mathcal{O}_{X'})$ lives in degree $\codim Y$ when $Y$ is a sub-bundle of $X'$ and in general we have the following diagram
\begin{center}
	\begin{tikzcd}
		H^{DR}_{-i}(0\times F/X) \arrow[r, "\sim"] \arrow[d] & {H^{2m+i}_{0\times F}(X\times F, \Omega^{\bullet}_{X\times F/X})} \arrow[d] \\
		H^{DR}_{-i}(X^u/X) \arrow[r, "\sim"] & {H^i_{X^u}(X, \mathcal{O}_X)}
	\end{tikzcd}
\end{center}
where $m$ is the relative dimension of $X \times F/X $ which is $\dim F$. At worst $F$ is the complete flag variety of $V$ which has dimension $\frac{k(k-1)}{2}$ and as $0\times F$ has codimension $nk$, we have a map
$$ H^{nk}_{0\times F}(X\times F, \Omega^*_{X \times F/X}) \rightarrow H^{nk - 2m}_{X^u}(X, \mathcal{O}_X)$$
and $m \leq \frac{k(k-1)}{2}$ so at worst this maps into degree $nk - k(k-1) = nk - k^2 + k > h$. We also have
$$\Omega^*_{X \times F / X} =  \mathcal{O}_{X \times F / X}\medoplus \Omega^1_{X \times F / X}[-1] \medoplus \cdots \medoplus \Omega^{\dim F}_{X \times F / X}[-\dim F]$$
and $H^i(Y, \mathcal{F}[-j]) = H^{i-j}(Y, \mathcal{F})$, so $H^{nk}_{0\times F}(X\times F, \Omega^j_{X \times F/X}[-j])$ appears in degree $nk + j$ and so gets mapped to $H^{nk - 2m+ j}_{X^u}(X, \mathcal{O}_X)$. This shows that we can ignore the local cohomology of any of the terms where $Y= 0 \times G/P$ where $P$ is any parabolic, including $P=G$.\par
Next we calculate what representations appear in the remaining terms, to do this we will use Borel-Weil-Bott and Lemma \ref{lem:localcoh_for_subbundles}. Consider the situation $X \times F \xrightarrow{q} F \xrightarrow{\pi} pt $. We have 
$$H^{\bullet}_Y(X \times F, \Omega^{\bullet}_{X\times F/X}) \cong H^{\bullet}(F, \mathcal{H}^{\bullet}_Y(X\times F, \mathcal{O}_{X \times F}) \medotimes \Omega^{\bullet}_F)$$
using $\Omega^{\bullet}_{X\times F/X} = \mathcal{O}_{X\times F} \medotimes q^*\Omega^*_F$, the projection formula and that calculating cohomology is the same pushing forward along $\pi \circ q$. Recall that we are using $\mathcal{H}^{\bullet}_Y(X\times F,-)$ to denote the pushdown along $q$ to $F$. We will first calculate what potential representations we get from the $\mathcal{H}^{\bullet}_Y(X\times F, \mathcal{O}_{X \times F})$ term, then we will add in $\Omega^*_F$ to get the result. 
\begin{lemma}\label{lem:localcoh_on_partial_flags}
	Let $Y = \Hom\left(W, V\right)_F$ where $F$ is some flag variety and $\dim W=k-l$. Then $\mathcal{H}^{\bullet}_Y(X \times F, \mathcal{O}_{X\times F})$ only contains terms that are the pushdown of line bundles $\mathcal{O}(\gamma)$ on the complete flag, where $\gamma = \sum a_i w_i$, with $a_{l} \geq n, a_i \geq 0$.
\end{lemma}
\begin{proof}
	We can push the bundle $\Hom\left(W, V\right)$ down to $Gr(S,W)$ without changing the weights that appear. (Another way of saying this is that the line bundle which pushes down to $\Hom\left(W, V\right)$ does not depend on the partial flag that it is over). So we can assume that $F = Gr(S,W)$. By Lemma \ref{lem:localcoh_for_subbundles} we have 
	$$\mathcal{H}^{\bullet}_Y(X\times Gr(S,W), \mathcal{O}_{X \times Gr(S,W)}) = \Sym{\bullet}\left(\Hom(W, V)^{*}\right)\medotimes \Sym{\bullet}\left(\Hom(L, V)\right)\medotimes \left(\det L^*\right)^{\otimes n}$$
	where $L = \ker (S \rightarrow W)$. Consider the first term, let $\alpha_{1} \geq \alpha_{2} \geq \cdots \geq \alpha_{k - l} \geq 0$, then using
	$$\Sym{\bullet}\left(\medoplus W\right) \cong \medotimes \Sym{\bullet}\left( W\right)$$ and the Littlewood-Richardson rule on 
	$$ \Sym{\alpha_{1}}\left( W\right) \medotimes \Sym{\alpha_{2}}\left( W\right) \medotimes \cdots \medotimes \Sym{\alpha_{k - l}}\left( W\right) \medotimes \Sym{0}\left( W\right) \medotimes \cdots \Sym{0}\left( W\right)$$
	we find that $\BBS^{\alpha}W$ appears for any $\alpha = (\alpha_{1}, \dots , \alpha_{l})$ a Young diagram. \par
	Looking at the final term, we have $\det L^* \cong \det W \cong \mathcal{O}(1)$. Now $\mathcal{O}(1)$ has Young diagram $(1,1, \dots , 1)$, therefore as we have $n$ copies of $\det L^*$, on the complete flag this corresponds to $nw_{l}$. Finally considering the middle term, we get $\BBS^{\beta} L^*$ appearing for any Young diagram $\beta$ in exactly the same way as for $W$.\par
	Putting these three pieces together we get 
	$$\BBS^{\beta}L^* \medotimes \BBS^{\alpha} W (n).$$
	Rewriting in terms of $W^*$, we get $\BBS^{\alpha'} W^*$ where $\alpha' = (-\alpha_{k-l}, \dots , -\alpha_1)$. This comes from the line bundle $\mathcal{O}(\gamma) = \mathcal{O}(\beta, \alpha') + \mathcal{0}(n,\dots , n, 0, \dots , 0)$ on the complete flag. We are allowed to shift by $(1,1, \dots, 1)$ as we only care about $SL(V)$ weights, not $GL(V)$ weights, so shift by $\alpha'_{1}$ to get that $\gamma$ must appear in the wanted form.
\end{proof}
\begin{proof}[Proof of Lemma \ref{lem:representationsappearinginlocalcohomology}]
	If a representation appears in $H_{X^u}^i ( X, \mathcal{O}_X)$, it must appear somewhere in the spectral sequence, i.e. it must appear in one of
	$$H^{\bullet}_Y(X \times F, \Omega^{\bullet}_{X\times F/X}) \cong H^{\bullet}(F, \mathcal{H}^{\bullet}_Y(X\times F, \mathcal{O}_{X \times F}) \medotimes \Omega^{\bullet}_F)$$
	for some $Y,F$ as above. By the earlier discussion we do not need to consider the terms where $Y = 0 \times F$ for any $F$. \par
	In Lemma \ref{lem:localcoh_on_partial_flags} we showed how to deal with $\mathcal{H}^{\bullet}_Y(X\times F, \mathcal{O}_{X \times F})$, so we next need to deal with $\Omega^{\bullet}_F$, where $\mathcal{O}_{F}\medoplus \Omega^1_{F}[-1] \medoplus \cdots \medoplus \Omega^p_{F}[-p]$, $\dim F = p$. In the case of $F=Gr(S,W)$ we have $\Omega^1_{Gr(S,W)} \cong \Hom(W, L) \cong W^* \medotimes L$ using the notation of Lemma \ref{lem:localcoh_on_partial_flags} and $\Omega^p_{Gr(S,W)} \cong \det (W^* \medotimes  L) \cong \mathcal{O}(-k)$. In the other extreme case, on the complete flag, $F(S)$, we have $\Omega^p_{F(S)} \cong \mathcal{O}(-2\rho)$ where $\rho = (k-1,k-2, \dots ,1,0)$. ($\rho$ is half the sum of the positive weights). \par
	In general there is not an easy way of describing $\Omega^{\bullet}_F$, but we can describe the line bundles that push down to the components of $\Omega^{\bullet}_F$. Using this description on $G/P$ we get that $\Omega^1_{G/P}$ has associated line bundles $\mathcal{O}(\gamma)$, where $\gamma$ corresponds to the collection of roots of $G$ with the roots of $P$ removed. As $B \subset P$ we get some collection of negative roots, and for $P=B$ we get all the negative roots. (The negative roots look like $(0,\dots ,0,-1,0, \dots ,0,1,0, \dots ,0)$)\\
	Taking wedge products corresponds to adding the weights together, for example on the complete flag, the canonical bundle corresponds to adding all the negative weights together, this is exactly $-2\rho$ as above. \par
	This means that we need to consider the line bundle corresponding to $nw_l + \sum a_i w_i + \sum_{p \in S^-} p$ where $S^-$ is some collection of negative roots. This is not written in the right form so we need to rewrite without the $\sum_{p \in S^-} p$ term. We can assume all the $a_i$ are zero as they only make the situation better, so let $\beta = n w_l + \sum_{p \in S^-} p$.\par
	Before we can analyse this further we will prove a bound on the maximum difference that can appear in $\gamma = \sum_{p \in S^-} p$. The claim is that in this situation we have $\gamma_j - \gamma_i \leq k + j - i - 1$, for $j > i$. We will show this by induction. \par
	Fix $k$, pick $i<j$. If $i=1, j=k$ then the weight $-2\rho$ gives the worst case. If at least one of $i,j$ is not $1,k$ then let $k' = j - i + 1 < k$. By assumption, only using weights that fit into the window $[i,j]$, we get that $\gamma'_j - \gamma'_i \leq k' + j - i - 1 = 2j - 2i$. We can at most increase $\gamma'_j$ by $i-1$ as we have that many choices for a negative weight with the $-1$ in one of the first $i-1$ entries and the $1$ in the $j^{th}$ entry, similarly we can reduce $\gamma'_i$ by at most $k-j$. Putting these together we get
	$$ \gamma_j - \gamma_i \leq \gamma_j' - \gamma'_i + i-1  + k - j \leq 2j - 2i + i + k - j - 1 = k + j - i - 1.$$
	The base case of $k=2$ is trivial.\par
	Now that we have this bound we can complete the proof. If $\beta$ is dominant, at most we can have $\beta_l - \beta_{l+1} = n - (k+l+1-l-1) = n-k$. Rewriting in terms of the $w_i$ basis we get the result.\par
	If $\beta$ is not dominant then we need to apply Borel-Weil-Bott in full generality. Let $\beta' = \sum_{\rho \in S^-} \rho$ and let $\beta_i$ be the smallest entry with $i \leq l, \beta_j$ the largest entry with $j > l$. After applying the twisted Weyl action from Borel-Weil-Bott we get $\alpha$ (we can assume that we get a dominant weight $\alpha$, if we get nothing it is fine) and we have the following bounds:
	\begin{align*}
	\alpha_{l} &\geq \beta_i + (k-i) - (k-l) = \beta_i + l - i,\\
	\alpha_{l+1} &\leq \beta_j +(k-j) -(k-(l+1)) = \beta_j + l - j + 1.
	\end{align*}
	As $\beta_i = n + \beta'_{i}$ and $\beta_j = \beta'_j$ we have
	\begin{align*}
	\alpha_l - \alpha_{l+1} &\geq n + \beta'_i + l - i - (\beta'_j + l - j - 1)\\
	&\geq n + (-k - j + i + 1) + l - i - l +j -1 \\
	&\geq n - k.
	\end{align*}
	Again, rewriting in the $w_i$ basis gives us the result.	
\end{proof}

\begin{remark}
	It is easy to show that these bounds are attained and that given any representation such that at least one of the coefficients in front of the $w_i$ is at least $n-k$ we can find a term in the spectral sequence where it appears. \cite[Theorem 2.1]{van1999local} enables us to strengthen this and say that these representations do appear in $H^i_{X^u}(X, \mathcal{O}_X)$, if a certain condition for the $SL(S)$ action holds. In these cases we get the same result as \cite[Theorem 4.6]{raicu2014local}. In particular, if $n > 3 \dim SL(S)$ then the condition holds, but it seems likely this bound can be strengthened, in particular, for $\dim S = 3$, it holds in all the cases we consider. (\cite{van1999local} does this situation as an example.) 
\end{remark}
\begin{remark}
	We can also see exactly where we get a stronger result then the one from \cite{van1991cohen}. In the above work we have shown that the collection of representations that we get does not get any bigger after applying Borel-Weil-Bott in full generality. In other words, everything that appears after applying the twisted Weyl group action can also be found directly. This is not true in \cite{van1991cohen}, more precisely, we have strengthened Corollary 6.8 of that paper.
\end{remark}
Now that we have found out what weights are Cohen-Macaulay, we can answer the original question of whether our specific choice of weights give a Cohen-Macaulay algebra.\\
To do this we use a result about what weights can appear in the decomposition of $\BBS^{\alpha} \medotimes \BBS^{\beta}$, the full decomposition is given by the Littlewood-Richardson rule, we just need the following.
\begin{lemma}\label{lem:littlewoodrichardsonestimates}
	Assume that $\BBS^{\gamma} \subset \BBS^{\alpha} \medotimes \BBS^{\beta}$ and that all diagrams have at most $k$ rows. Then we have the following bounds on the entries of $\gamma$,
	$$ \alpha_i + \beta_k \leq \gamma_i \leq \alpha_1 + \beta_i.$$
\end{lemma}
\begin{proof}
	We will first prove the upper bound then the lower one.	A straightforward application of the rule says that the boxes from the $i^{th}$ row of $\beta$ can only appear in rows $i$ to $k$. This gives us $\gamma_{1} \leq \alpha_{1} + \beta_1$. As each column has to have strictly increasing boxes\footnote{Part of the L-R rule says you need to label the boxes of $\beta$ by the row that they appear in. Then the columns of $\gamma$ need to be strictly increasing.}, once we have filled up to column $\alpha_{1}$ the columns have to look like $1 2 3 4$ etc. This shows that row $i$ can not be longer than $\alpha_{1} + \beta_i$. \par
	Using the first observation also tells us that $\gamma_{k} \geq \alpha_{k}+ \beta_{k}$ and we can generalize this using part of the rule. It tells us that there must be at least $\beta_{k}$ boxes from $\beta_{k-1}$ placed in the rows above the $k^{th}$ row, combining this with the fact that boxes from $\beta_{k-1}$ can only appear in rows $k-1$ and $k$ tells us that at least $\beta_{k}$ of them must appear in row $k-1$, therefore $\gamma_{k-1} \geq \alpha_{k-1} + \beta_{k}$. In exactly the same way one can show that we need to place at least $\beta_{k}$ boxes from $\beta_{i}$ in the $i^{th}$ row. This gives us the lower bound $\gamma_{i} \geq \alpha_{i} + \beta_{k}$. 
\end{proof}
This is enough to show that our algebra is Cohen-Macaulay.
\begin{proposition}\label{prop:algebraisCM}
	The algebra
	$$\Lambda = \End \left(\bigoplus_{\alpha \in \mathcal{UP}_{n,k}} (\BBS^{\alpha}S^{*}\medotimes \Sym{\bullet}X^*)^G \right)$$
	is Cohen-Macaulay.
\end{proposition}
\begin{proof}
	It is sufficient to show that every piece of the decomposition of $\Lambda$ is Cohen-Macaulay. In this situation each piece appears in the decomposition of $\BBS^{\alpha^*} S^*\medotimes \BBS^{\beta}S^*$ where $\alpha^* = (\alpha_1 - \alpha_k, \alpha_1 -\alpha_{k-1}, \dots , \alpha_1 - \alpha_1)$ and $\alpha_i, \beta_i \leq (n-k)\frac{k-i}{k}$. Applying Lemma \ref{lem:littlewoodrichardsonestimates} we get the following bounds on the $\gamma_i$, 
	$$ \alpha_1 - \alpha_{k-i +1}\leq \gamma_i \leq \alpha_1 + \beta_i .$$
	(as $\alpha_k = \beta_k = 0$). This shows that 
	\begin{align*}
	\gamma_i - \gamma_{i+1} &\leq \alpha_1 + \beta_i  - (\alpha_1 - \alpha_{k-i}) \\
	&\leq \beta_i + \alpha_{k-i} \\
	&\leq (n-k)\frac{k-i}{k} + (n-k) \frac{k-(k-i)}{k} = n-k.
	\end{align*} 
	As $n$ and $k$ are coprime we have a strict less than and therefore if we write $\gamma$ in terms of the $w_i$ each coefficient must be less then $n-k$. Therefore our algebra, $\Lambda$, is Cohen-Macaulay. 
\end{proof}
\noindent
As an aside, in general a collection of weights will give a Cohen-Macaulay endomorphism algebra if it is "small" enough and that algebra will have finite global dimension if the collection is "big" enough. As we are looking for a NCCR which requires both Cohen-Macaulay and finite global dimension, our collection should be on the boundary of both conditions. This was idea can be formalised by maximal modification algebras, \cite{2010arXiv1007.1296I}. In particular they prove that NCCR's are maximal modifying algebras, \cite[Prop. 4.5]{2010arXiv1007.1296I}, this means that adding any new weight to our collection will stop our algebra being Cohen-Macaulay.\par
It turns out that we can also prove this directly in our case, without assuming finite global dimension. To prove this, we first need another simple result about the Littlewood-Richardson rule, it follows by simply building $\beta$ given $\alpha$ and $\gamma$.
\begin{lemma} \label{lem:maxdualdiagram}
	Let $\alpha, \beta, \gamma$ be weights with at most $l$ rows, then a $\BBS^{\gamma}$ such that $\gamma_1 = \alpha_1$ appears in the decomposition of $\BBS^{\alpha} \medotimes \BBS^{\beta}$ if and only if $\beta \leq (\alpha_1 - \alpha_l, \alpha_1 - \alpha_{l-1}, \dots , \alpha_1 - \alpha_1)$.
\end{lemma}
\begin{lemma}
	Let $\beta$ be any weight that does not appear in $\mathcal{UP}_{n,k}$, i.e. there exists a $\beta_i$ that does not satisfy $\beta_i \leq (n-k)\frac{k-i}{k}.$ Then there exists a weight $\alpha \in \mathcal{UP}_{n,k}$ such that at least one representation that appears in the decomposition of $\BBS^{\alpha^*} \medotimes \BBS^{\beta}$ is not Cohen-Macaulay.
\end{lemma}
\begin{proof}
	Let $l = \max \{ i | \beta_i > (n-k)\frac{k-i}{k}\}$ by the assumptions on $\beta$ such a $l$ exists. Let 
	$$\bar{\alpha} = (\alpha_1, \alpha_2, \cdots \alpha_{l}, 0, \dots , 0)$$
	where $\alpha_i$ is the largest allowed value. Then we have that 
	$$\bar{\alpha}^* = (\alpha_1, \alpha_1, \dots , \alpha_1, \alpha_1 - \alpha_{k-l}, \alpha_1 - \alpha_{k-l -1}, \dots , \alpha_1 - \alpha_2, 0).$$
	We claim that we can find $\BBS^{\gamma} \subset \BBS^{\bar{\alpha}^*} \medotimes \BBS^{\beta}$ such that $\BBS^{\gamma}$ is not Cohen-Macaulay. It is sufficient to show that there exists a $\gamma$ with $\gamma_l - \gamma_{l+1} \geq n-k$. \par
	Let $\gamma = (\alpha_1 + \beta_1, \alpha_1 + \beta_2, \dots , \alpha_1 + \beta_l, ?, \dots , ?)$. In applying the Littlewood Richardson rule (algorithm) we have used up all the entries of $\beta$ up to and including the $l^{th}$ row, we also can not change the first $l$ entries of $\gamma$ so we have reduced the problem to applying the Littlewood Richardson rule to the situation 
	$$(\alpha_1 - \alpha_{k-l}, \alpha_1 - \alpha_{k-l -1}, \dots , \alpha_1 - \alpha_2, 0)\otimes (\beta_{l+1}, \dots , \beta_{k-1}, 0).$$
	We want to find a diagram in the decomposition with the first entry equal to $\alpha_1 - \alpha_{k-l}$. To check that we can do this we need to show that our weights satisfy the conditions of the above lemma. \par
	So we need to show that $\beta_{l+i} \leq (\alpha_1 - \alpha_{k-l}) - (\alpha_1 - \alpha_{i}) = \alpha_{i} - \alpha_{k-l}$. Note that we have $\alpha_i + \alpha_{k-j} \leq \alpha_{i-j}$ for $i \geq j$. (We are assuming that $\alpha_{i}$ is the maximum possible value for all $i$.)\footnote{This follows as $\lfloor a \rfloor + b \leq a+b$ which implies that $\lfloor a \rfloor + \lfloor b \rfloor = \lfloor \lfloor a \rfloor + b \rfloor \leq \lfloor a+b \rfloor$.} Applying this gives us that 
	$$\alpha_{i} - \alpha_{k-l} \geq \alpha_{l+i}$$
	and by assumption $$\beta_{l+i} \leq (n-k)\frac{k-(l+i)}{k}$$
	so using that $\alpha_{l+i}$ is as large as possible, we get that $\beta_{l+i} \leq \alpha_{l+i}$ and so we have found a $\gamma$ which looks like $(\dots, \alpha_{1} + \beta_l, \alpha_{1} - \alpha_{k-l}, \dots)$. And we have $\beta_l + \alpha_{k-1} \geq n-k$, use that $\alpha_{k-l}$ is maximal and that $\beta_l > (n-k)\frac{k-l}{k}$.
\end{proof}
\subsection{Finite global dimension}\label{sec:fgdim}
To complete the proof that $\Lambda$ is a NCCR, we need to show that $\Lambda$ has finite global dimension as a graded algebra. To do this we will use the resolution from \cite{2015arXiv150205240S} described at the start of this section. The same resolution is also found in \cite{fonarev2013minimal} and \cite{donovan2014window} for the special case of $G=GL(S)$ but they describe it differently (called staircase complexes in those papers). It turns out that their description is more helpful for us.\par
We will analyse this resolution and show that by iterating it a finite number of times, any module $P_{\alpha}$ has a finite resolution by modules $P_{\beta}$ for $\beta \in \mathcal{UP}_{n,k}$. First if $P_{\alpha}$ is any module such that $\alpha_{1} > n-k$, then using the description in \cite{2015arXiv150205240S} it is easy to see that all the $P_{\beta}$ that appear in the resolution satisfy $\beta_{1} < \alpha_{1}$, see Theorem \ref{thm:FGDandNCCR} for the proof.\\
This tells us that every weight has a finite resolution by weights that fit into a rectangle of size $k\times (n-k)$. Therefore it is sufficient to prove that all weights of size at most $k\times (n-k)$ have finite projective dimension. \par
To do this we will use the description of the resolution, as found in \cite{fonarev2013minimal} and \cite{donovan2014window}, called a staircase complex. This is all found in \cite{fonarev2013minimal}, we reproduce it here for completeness, and also simplify the argument. Fonarev uses $GL(S)$ representations, our representations are $SL(S)$ ones.\par
We can represent a Young diagram, $\alpha$, of size $k\times (n-k)$ using a binary sequence of length $n$ containing $k$ 1's and $n-k$ 0's. This is done by considering the path from the bottom left corner to the top right corner of the Young diagram, let a vertical section be a 1 and a horizontal section be a 0. As an example, consider $\alpha = (3,1,0)$, thought of as a Young diagram of size $3 \times 4$, the binary sequence is 1010010, another example can be found after Figure 2. Note that every binary sequence of length $n$ with $k$ 1's and $n-k$ 0's correspondences uniquely to a Young diagram of size at most $k \times (n-k)$. Representing diagrams in this way gives us a natural action of $\BBZ/n\BBZ$ by rotation and under this action exactly one diagram in each orbit will be upper triangular. \par
There is a geometric way of describing this action as well. Take a diagram $\alpha$ with $\alpha_{1} \leq n-k$ and glue them together at the corners. (In terms of the binary sequence, repeat the sequence infinitely in both directions.) Then applying $i$ rotations is the same as starting at the bottom left corner of the diagram and taking $i$ steps along the edge of $\alpha$, then creating a new diagram with the bottom left corner being at this point. I.e. take a box of size $k\times (n-k)$ and move it along the edge of the diagram $\alpha$ to give the other diagrams in the orbit. \par
\begin{figure}[h!]
	\centering
	\includegraphics[width=.8\linewidth]{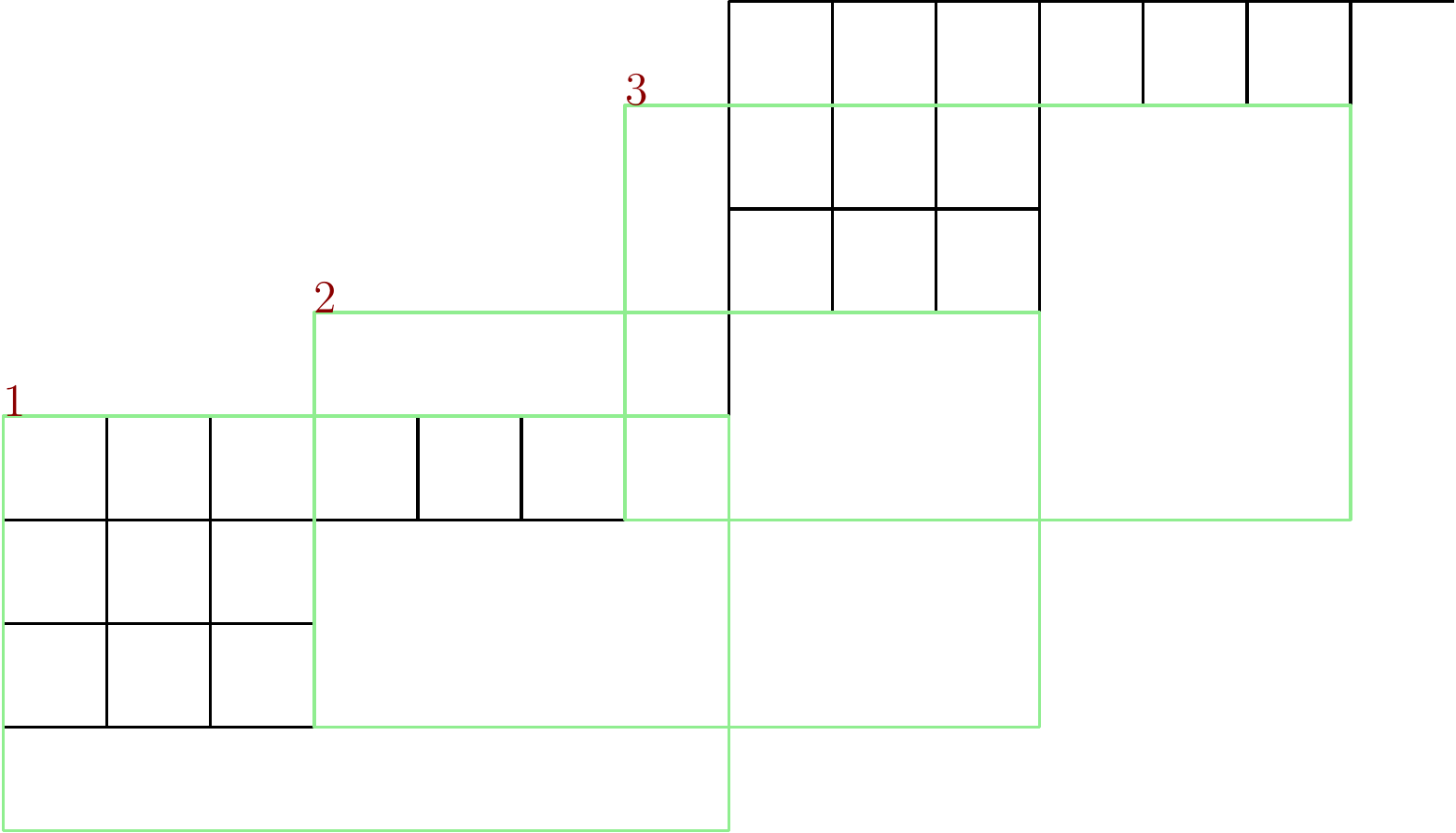}
	\caption{The effect of the $\BBZ/n\BBZ$ action}
\end{figure}
\noindent
In Figure 2 we see an example where $n=11, k=4, \alpha =(6,3,3,0)$ and the binary sequence is 10001100010. The green box labelled 1 is the original diagram, the box labelled 2 is the diagram $(4,3,0,0)$ and corresponds to $4 \cdot \alpha$. Finally, the box labelled 3 is the diagram $(4,4,1,0)$ and corresponds to $9 \cdot \alpha$. \par
Using this description we can also easily show that there exists a unique upper triangular diagram in the orbit. To see this draw a diagonal line with gradient $(n-k)/k$ through each corner of the diagram, one of these lines will be the lowest. In fact only one, using the fact that $n$ and $k$ are coprime. Take the diagram that has as corners, the points where this line touches the repeating diagram, it will be upper triangular. See Figure 3 for an example, again with $\alpha = (6,3,3,0)$.
\begin{figure}[h!]
	\centering
	\includegraphics[width=.8\linewidth]{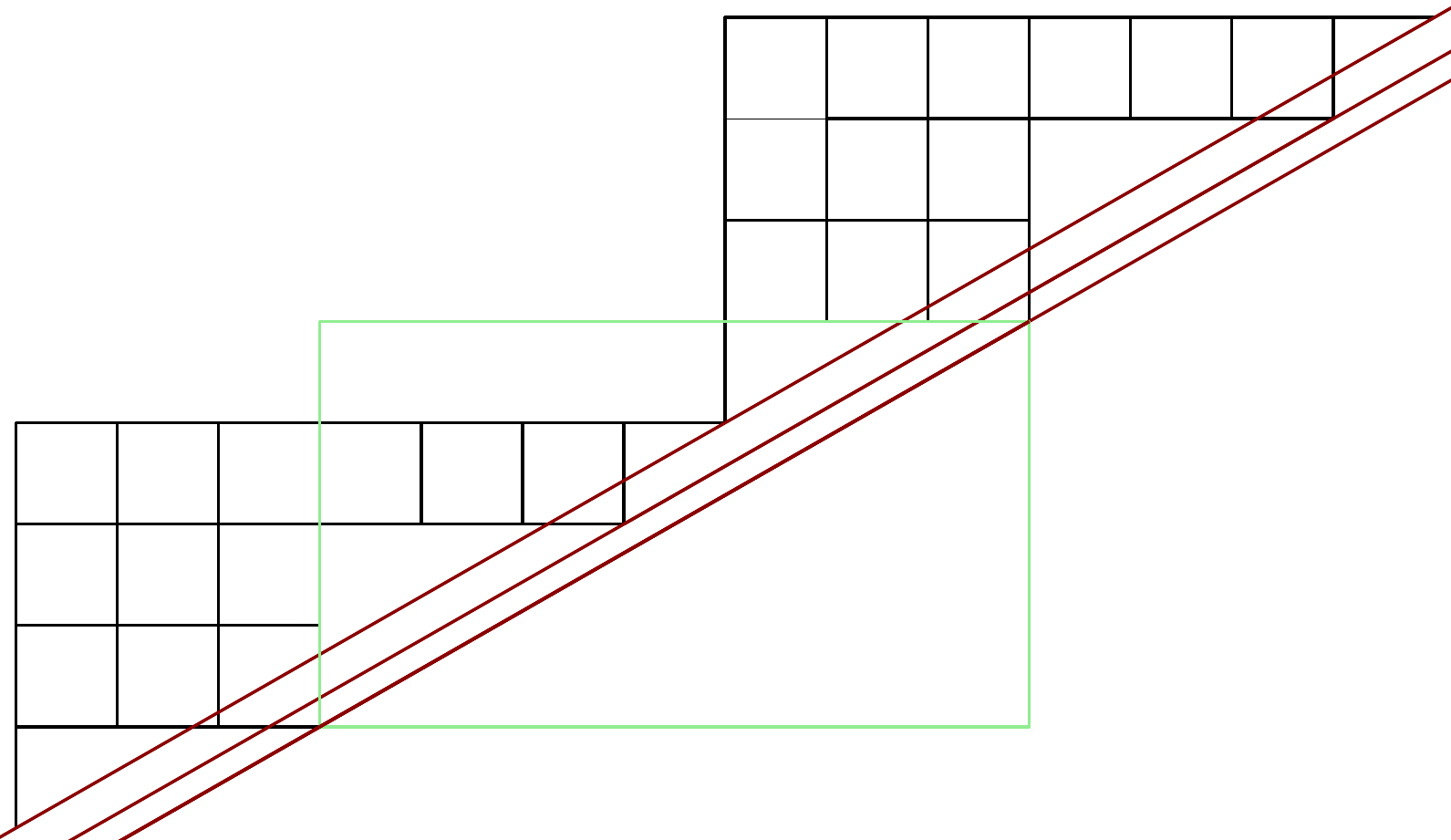}
	\caption{Detecting the upper triangular element in the orbit}
\end{figure}
Let $d^{upp}(\alpha)$ be the number of rotations needed to take $\alpha$ to an upper triangular diagram, or equivalently the number of steps needed to move to the corner which gives the lowest line. Note that $d^{upp}(\alpha) = 0$ for all diagrams $\alpha \in \mathcal{UP}_{n,k}$.\par
Next we describe the terms that appear in the staircase complex, these come from $SL(S)$ representations, but we will not remove any complete columns when we draw the diagrams. \\
Let $\alpha$ be a diagram such that $\alpha_{1} \leq n-k$, create a new diagram by adding a strip to the edge of $\alpha$ as seen in Figure 4 and extend along the first row until we have added $n$ boxes. Then the first term in the complex has diagram corresponding to the new diagram, but where we have only completed the first column, and the $i^{th}$ term has diagram corresponding to the new diagram, but only the new terms from the first $i$ columns. It is also possible to find the multiplicities of the terms, but we do not need that level of detail. Looking at Figure 4, the $i^{th}$ term is the diagram corresponding to the original diagram as well as all boxes labelled $i$ or less.
\begin{figure}[h!]
	\centering
	\includegraphics[width=.6\linewidth]{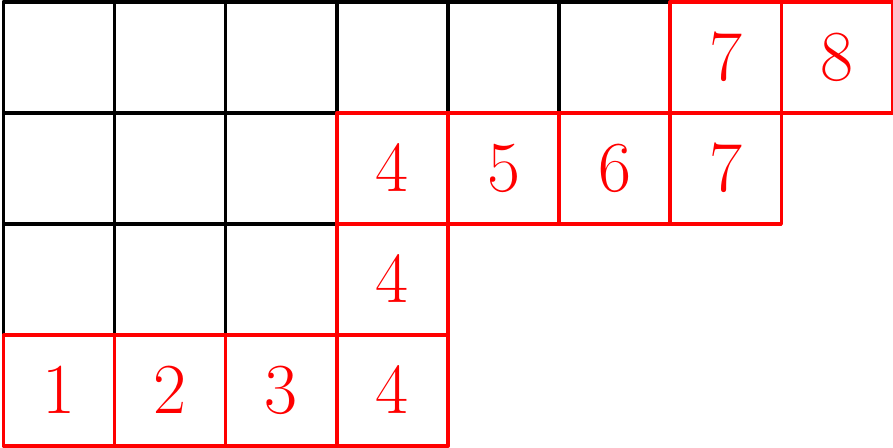}
	\caption{The staircase resolution}
\end{figure}
\newpage
\begin{lemma}\label{lem:metricforresolution}
	Let $\alpha$ be a $SL_k$ representation whose Young diagram is of size $k \times (n-k)$ and not an element of $\mathcal{UP}_{n,k}$. Let $P_{\beta}$ be any term that appears in the staircase resolution of $P_{\alpha}$.	Then we have $d^{upp}(\beta)< d^{upp}(\alpha)$.
\end{lemma}
\begin{proof}
	We will take the definition of $d^{upp}$ as the number of steps to get to the corner that gives the lowest diagonal. It is clear that adding more boxes can only make the lowest diagonal lower so we can assume $\beta$ is the final term in the resolution. We are working with $SL_k$ representations, so we can remove any complete columns from a Young diagram, and if we are calculating $d^{upp}(\gamma)$ for $P_{\gamma}$ then $\gamma$ has no complete columns.\par
	Let $\overline{\beta}$ be the final Young diagram in the staircase resolution before we remove any columns. We have $d^{upp}(\beta) \leq d^{upp}(\overline{\beta})$ as removing columns reduces $d^{upp}$, in fact by construction, $\overline{\beta}$ has a complete column, so let $\beta'$ be the Young diagram corresponding to $\overline{\beta}$ with the first column removed, we then still have $d^{upp}(\beta) \leq d^{upp}(\beta')$. \par
	Therefore it is sufficient to show $d^{upp}(\beta') < d^{upp}(\alpha)$. We have $\beta' = (a, \alpha_{1}, \dots , \alpha_{k-1})$ for some $a \geq \alpha_{1}$, so if the lowest diagonal comes from a corner not on the first row then $d^{upp}(\beta') = d^{upp}(\alpha) - 1$, as we need to make one less step to reach the lowest diagonal on $\beta'$ then on $\alpha$ (recall, $\alpha_{k}=0$ by definition).\par
	If the lowest diagonal for $\beta'$ touches a corner on the first row, we get $(\alpha_{1} + i, \dots , \alpha_{k-1}+i, i)$ is upper triangular, which implies $d^{upp}(\alpha) =0$. By assumption this is not the case and therefore
	$$ d^{upp}(\beta) \leq d^{upp}(\beta') = d^{upp}(\alpha) -1 < d^{upp}(\alpha).$$
\end{proof}
\noindent
This result is found in \cite[Prop. 5.6]{fonarev2013minimal} as a note. We can now prove that $\Lambda$ is a NCCR for $X/G$. 
\begin{theorem}\label{thm:FGDandNCCR}
	The graded algebra 
	$$\Lambda = \End \left(\bigoplus_{\alpha \in \mathcal{UP}_{n,k}} (\BBS^{\alpha}S^{*}\medotimes \Sym{\bullet}X^*)^G \right)$$
	has finite global dimension and is a NCCR of $X/G$.
\end{theorem}
\begin{proof}	
	It is sufficient to show that $P_{\alpha}$ has finite projective dimension for any Young diagram $\alpha$. (Basically by graded Nakayama's lemma, as the $P_{\alpha}$ resolve the graded simples.) By construction if $\alpha \in \mathcal{UP}_{n,k}$ then $P_{\alpha}$ is projective. Therefore let $\alpha$ be a weight not in $\mathcal{UP}_{n,k}$. If $\alpha_{1}\leq n-k$ then Lemma \ref{lem:metricforresolution} shows that we have a finite resolution of $P_{\alpha}$ by projectives. \par
	Another way of describing the terms, $P_{\beta}$, in the resolution is that $\beta$ is one of the weights that arise by applying Borel-Weil-Bott to $\alpha + (0,0, \dots , 0 , i)$ for $i=1, \dots n$. We will use this description to deal with diagrams $\alpha$ such that $\alpha_{1} > n-k$.\par
	As $i \leq n < \alpha_1 + k$, after applying Borel-Weil-Bott we still have that $\alpha_{1}$ is the largest entry, and the final element is at least 1. This shows that any $P_{\beta}$ in the resolution satisfies $\beta_1 < \alpha_{1}$. Therefore any module $P_{\alpha}$ with $\alpha_1 > n-k$ has a finite resolution by modules $P_{\beta}$ such that $\beta_{1} \leq n-k$.\par
	Putting these two results together we get that each $P_{\alpha}$ has finite projective dimension and therefore $\Lambda$ has finite global dimension. We have already shown that $\Lambda$ is Cohen-Macaulay in Proposition \ref{prop:algebraisCM} so we get $\Lambda$ is an NCCR of $X/G$.
\end{proof}

\section{Equivalences between NCCR's}\label{sec:equivNCCR}
Observing that the Grassmannian can be equivalently described as either the space of subspaces or the space of quotients gives us another NCCR. This observation is also the motivation for Hori duality for the group $SL$.\par
Let $V,S$ and $Q$ be vector spaces with dimension $n, k$ and $n-k$, then $Gr(S, V) = Gr(V, Q)$ and $\Hom(S,V)/SL(S) \cong \Hom(V,Q)/SL(Q)$. A straightforward modification of Section \ref{sec:algebraicnccr} proves that 
$$\Lambda' = \left(\BBS^{\alpha^t} Q^* \medotimes \Sym{\bullet} X'^*\right)^{SL(Q)}$$ where $X' = \Hom(V,Q)$ and the sum is over exactly the same Young diagrams, $\alpha \in \mathcal{UP}_{n,k}$ is an NCCR of $\Hom(V,Q)/SL(Q)$.\par
As these spaces are isomorphic we actually get another NCCR of $\Hom(S,V)/SL(S)$, however we can show that even in the first non-trivial case ($n=5, k=2$) we get two different endomorphism algebras. Representing the algebras as quivers and writing $X=\Hom(S,V)$ and $X' = \Hom(V,Q)$ we get
\begin{center}
	\begin{tikzcd}
		\cdot  \arrow[loop, out = -150, in =150,distance = 4em,"{\medoplus \BBS^{(a, a)} V^*}"] \arrow[rr, "{\medoplus \BBS^{(a+1, a)} V^*}"', bend right] &  & \cdot \arrow[ll, "{\medoplus \BBS^{(a+1, a)} V^*}"', bend right] \arrow[loop, out = 30, in =-30,distance = 4em, "{\medoplus \BBS^{(a, a)} V^* \medoplus \BBS^{(a+2, a)} V^*}"]
	\end{tikzcd}
\end{center}
for $X/SL(S)$ (we have not mentioned what the relations are) and
\begin{center}
	\begin{tikzcd}
		\cdot  \arrow[loop, out = -150, in =150,distance = 4em,"{\medoplus \BBS^{(a, a, a)} V}"] \arrow[rr, "{\medoplus \BBS^{(a+1, a+1,a)} V}"', bend right] &  & \cdot \arrow[ll, "{\medoplus \BBS^{(a+1, a, a)} V}"', bend right] \arrow[loop, out = 30, in =-30,distance = 4em, "{\medoplus \BBS^{(a, a, a)} V \medoplus \BBS^{(a+2, a+1 ,a)} V}"]
	\end{tikzcd}
\end{center}
for $X'/SL(Q)$. Picking a volume form for $V$ gives us an isomorphism between some of the components ($\medoplus \BBS^{(a, a, a)} V \cong \medoplus \BBS^{(a, a)} V^*$), but not all of them. For example, consider the $a=0$ component, we have $|V^*|$ maps from the left to the right in the first algebra and $|\medwedge^2 V|$ maps in the second algebra. These have ranks 5 and 10 respectively. \par
We want to prove that these algebras are derived equivalent. This is expected for at least two reasons. First, all NCCR's are conjectured to be derived equivalent, in fact all crepant categorical resolutions are conjectured to be equivalent and second, it is predicted by Hori duality.
\subsection{Geometric Resolution}
In the process of showing the derived equivalence between the algebras we will introduce a crepant categorical resolution and prove that it is equivalent to both algebras. This equivalence will later turn out to be a key part of the HPD result, see Section \ref{sec:HPDfibre}. This categorical resolution uses the methods of \cite{2006math......9240K}, the idea is to first take a geometric resolution of $X/G$. In our case, it is given by $\mathcal{O}_{Gr(S,V)}(-1)$ with exceptional divisor, the zero locus, isomorphic to the Grassmannian. Let $Gr = Gr(S,V)$. We have
\begin{center}
	\begin{tikzcd}
		{Gr} \arrow[d] \arrow[r, "i", hook, shift left] & {\mathcal{O}_{Gr}(-1)} \arrow[d, "\pi"] \arrow[l, "p", shift left] \\
		0 \arrow[r, hook]                        & {\Hom(S, V)/SL(S)}                                    
	\end{tikzcd}
\end{center}
and $D^b(\mathcal{O}_{Gr}(-1))$ is a categorical resolution. We also have a grading as we have a $\BBC^*$-action which is scaling the fibres, this action is a lift of the $\BBC^*$ action on the singularity. \par
We want to find a smaller subcategory that is still a categorical resolution and is in addition crepant. To do this we use the rectangular Lefschetz decomposition of $D^b(Gr)$, found by Fonarev, \cite{fonarev2013minimal}.
$$ D^b(Gr) = \langle\mathcal{A}, \mathcal{A}(1), \dots , \mathcal{A}(n-1) \rangle $$
The first block, $\mathcal{A}$, is given by Schur powers of the universal sub-bundle, $S^*$, on the Grassmannian, indexed by upper triangular Young diagrams of size at most $(n-k) \times k$, exactly the ones in $\mathcal{UP}_{n,k}$ as earlier. Given this, Kuznetsov \cite{2006math......9240K}, finds a partial decomposition of $D^b(\mathcal{O}_{Gr}(-1))$ and the left orthogonal of this collection gives a categorical (crepant) resolution.\par
Applying this in our situation gives us
$$\mathcal{B} = \{ \mathcal{F} | i^*\mathcal{F} \in \mathcal{A}\},$$
a categorical weakly crepant resolution\footnote{It is always a categorical resolution, it is crepant only when $n$ and $k$ are coprime.}. If in addition $\mathcal{A}$ is generated by the pullback of a vector bundle from $\mathcal{B}$, and the associated endomorphism algebra pushes down along $\pi$ to a sheaf of algebras, Kuznetsov proves that the pushdown of the endomorphism algebra gives an NC(C)R. Now as $p \circ i = id$, it is clear that $i^*\left(\medoplus p^* \BBS^{\alpha} S^*\right)$ generates $\mathcal{A}$. Therefore it is sufficient to show this bundle is tilting over $\Hom(S,V)/SL(S)$ to get another NCCR.
\begin{lemma}\label{lem:tiltingbundle}
	The vector bundle $\medoplus p^* \BBS^{\alpha}S^*$ where the sum is over $\alpha \in \mathcal{UP}_{n,k}$ is tilting over \\
	$\Hom(S, V)/SL(S)$.
\end{lemma}
This Lemma actually holds for all $n,k$, where we allow non-strictly upper triangular elements.
\begin{proof}
	As $\Hom(S, V)/SL(S)$ is affine we just need to show that 
	$$R\Hom(p^* \BBS^{\alpha}S^*, p^* \BBS^{\beta}S^*)= \Hom(p^* \BBS^{\alpha}S^*, p^* \BBS^{\beta}S^*)$$
	for all $\alpha, \beta \in \mathcal{UP}_{n,k}$. We have
	\begin{align*}
		R\Hom(p^* \BBS^{\alpha}S^*, p^* \BBS^{\beta}S^*) &\cong H^{\bullet}\left(\mathcal{O}(-1), p^*\left(\BBS^{\alpha}S \medotimes \BBS^{\beta}S^*\right)\right)\\
		& \cong H^{\bullet}\left(Gr, \BBS^{\alpha}S \medotimes \BBS^{\beta}S^*\medotimes p_* \mathcal{O}_{\mathcal{O}(-1)}\right).
	\end{align*}
	As $p_* \mathcal{O}_{\mathcal{O}(-1)} = \medoplus_{i\geq 0} \mathcal{O}(i)$ we need to show that $\BBS^{\alpha}S \medotimes \BBS^{\beta}S^* \medotimes \mathcal{O}(i)$ has no higher cohomology for all $i\geq 0$. To do this consider what terms can appear in the decomposition of $$\BBS^{\alpha}S \medotimes \BBS^{\beta}S^* \medotimes \mathcal{O}(i) \cong \BBS^{\alpha^*}S^* \medotimes \BBS^{\beta}S^*\medotimes \mathcal{O}(i).$$
	From Lemma \ref{lem:littlewoodrichardsonestimates} and the bounds on $\alpha$ and $\beta$, $(\beta_k = 0)$ we get the following bounds on any $\BBS^{\gamma}S^*$ appearing
	$$ -(n-k)(k-1)/k +i \leq -\alpha_{1} + i  \leq \gamma_{k} \leq i.$$
	If $\gamma_{k}\geq 0$, it is a dominant weight and we only have cohomology in degree 0 (using Borel-Weil-Bott). Else we have 
	$$0 > \gamma_{k} \geq i - (n-k)(k-1)/k \geq - (n-k)(k-1)/k,$$ so when we apply the twisted Weyl action, we first add $\rho$ which has $k^{th}$ entry $n-k+1$, then we rearrange and as $n-k + 1 > (n-k)(k-1)/k$ we get two entries with the same value. Therefore this weight is not dominant and we have no cohomology at all. In both cases we get the wanted result.
\end{proof}
\begin{corollary}\label{cor:nccrexistence}
	$\pi_*\mathcal{E}nd \left(\medoplus p^* \BBS^{\alpha}S^*\right)$ is an NCCR for $\Hom(S,V)/SL(S)$.
\end{corollary}
\begin{proof}
	See \cite[Section 5]{2006math......9240K}.
\end{proof}
\begin{lemma}\label{lem:generatorsforcategoricalresolution}
	$\mathcal{B} = \langle p^* \BBS^{\alpha}S^*\rangle_{\alpha \in \mathcal{UP}_{n,k}}$
\end{lemma}
\begin{proof}
	This follows from Corollary \ref{cor:nccrexistence}. Kuznetsov proves $\pi_*\mathcal{E}nd \left(\medoplus p^* \BBS^{\alpha}S^*\right)$ is a NCCR by showing 
	$$\mathcal{B} \simeq D^b(\Hom(S,V)/SL(S), \pi_*\mathcal{E}nd \left(\medoplus p^* \BBS^{\alpha}S^*\right)).$$
	Now as $\Hom(S,V)/SL(S)$ is affine, we get $\mathcal{B} \simeq D^b(\End \left(\medoplus p^* \BBS^{\alpha}S^*\right))$ and therefore the wanted result follows. (Use the grading to show that the $p^* \BBS^{\alpha}S^*$ generated $\mathcal{B}$.)
\end{proof}
This process actually gives us four potentially different categorical resolutions. We have two different Lefschetz decompositions, one with terms of the form $\BBS^{\alpha}S^*$, and one with terms $\BBS^{\alpha^t}Q^*$, where $\alpha^t$ is the transposed Young diagram. We then get two more by using $\BBS^{\alpha}S$ and $\BBS^{\alpha^t}Q$. As we have $\left(\BBS^{\alpha}Q\right)^{*} = \BBS^{\alpha}Q^{*}$ and dualizing is an anti-auto-equivalence we only have two potentially different categorical resolutions. We can in fact go further and show that the remaining two categorical resolutions are the same.
\begin{lemma}\label{lem:differentgeneratingsets}
	The two collections $\{\BBS^{\alpha}S^*\}_{\alpha \in \mathcal{UP}_{n,k}}$ and $\{\BBS^{\alpha^t}Q^*\}_{\alpha \in \mathcal{UP}_{n,k}}$ generate the same subcategory in $D^b(Gr(S,V))$. 
\end{lemma}
\begin{proof}
	We can take Schur powers of a sequence, as follows. Given a short exact sequence
	$$0 \rightarrow A \rightarrow B \rightarrow C \rightarrow 0$$
	we get 
	\begin{align*}
	0 \rightarrow \BBS^{\alpha} A \rightarrow \BBS^{\alpha} B \rightarrow \cdots \rightarrow &\medoplus \left(\BBS^{\beta} B \medotimes \BBS^{\gamma^t} C\right)^{\oplus m^{\alpha}_{\beta \gamma}} \rightarrow \cdots \rightarrow \BBS^{\alpha^t} C \rightarrow 0
	\end{align*}
	is exact, where the sum is over $|\gamma| = p, |\beta| = |\alpha| - p$ and $m^{\alpha}_{\beta \gamma}$ is the Littlewood-Richardson coefficients. It is non-zero only if $\gamma$ is a sub-digram of $\alpha$. Applying this to
	$$0 \rightarrow Q^* \rightarrow V^* \medotimes \mathcal{O}_{Gr(S,V)} \rightarrow S^* \rightarrow 0$$
	we get that 
	$$\langle \BBS^{\alpha^t} Q^* | \alpha \in \mathcal{UP}_{n,k} \rangle \subset \langle \BBS^{\alpha} S^* |\alpha \in \mathcal{UP}_{n,k} \rangle.$$
	One can also resolve $\BBS^{\alpha}C$ by terms involving Schur powers of $A$ and $B$, this gives us the reverse inclusion. Therefore our two collections are just different generators of the same subcategory.
\end{proof}
\noindent
This shows that all four categorical resolutions that arise from the geometric resolution are equivalent, and therefore all four associated NCCR's are derived equivalent.
\subsection{Equivalence}
To prove that the algebraic NCCR is equivalent to the categorical crepant resolution we will use a larger stack which has open sets associated to both the algebraic and geometric NCCR's. We will then use a VGIT style argument to show that these two NCCR's are derived equivalent. \par
Consider the stack $\mathcal{Z} = \left[\Hom(S,V)\times \wedge^k S  /GL(S) \right]$. We find two opens inside $\mathcal{Z}$, that are actually GIT quotients for $\mathcal{Z}$. One is 
$$Y = \left[\Hom(S,V)^f\times \wedge^k S /GL(S)\right] \cong \mathcal{O}_{Gr}(-1)$$
where $\Hom(S,V)^f$ is all the full rank maps and the other is 
$$\mathcal{X} = \left[(\Hom(S,V)\times \wedge^k (S\setminus 0))/GL(S) \right] \cong \left[ \Hom(S,V)/SL(S)\right].$$
They fit together as follows. 
\begin{center}
	\[\begin{tikzcd}
	& {\mathcal{Z} =\left[\frac{\Hom(S,V)\times \wedge^k S}{GL(S)} \right]} \\
	{Y =\mathcal{O}_{Gr}(-1)} && {\mathcal{X} = \left[\frac{\Hom(S,V)}{SL(S)}\right]}
	\arrow["i", hook, from=2-1, to=1-2]
	\arrow["j"', hook', from=2-3, to=1-2]
	\end{tikzcd}\]
\end{center}
We will define a subcategory of $D^b(\mathcal{Z})$, that restricts to both the geometric NCCR, $\mathcal{B}$, along $i$ and to the algebraic NCCR, $D^b(\Lambda)$, along $j$. We will also prove that these restrictions are fully faithful on this subcategory, completing the proof. Let $\mathcal{C} \subset D^b(\mathcal{Z})$ be defined as, 
$$\mathcal{C} = \{\mathcal{O}_{\mathcal{Z}} \medotimes \BBS^{\alpha} S^* | \alpha \in \mathcal{UP}_{n,k} \}. $$
\begin{proposition}\label{prop:linebundleequivfibre}
	The subcategory $\mathcal{C}$ is equivalent to $\mathcal{B}$.
\end{proposition}
\begin{proof}
	$\mathcal{C}$ is defined by a set of generators and Lemma \ref{lem:generatorsforcategoricalresolution} gives us a set of generators for $\mathcal{B}$. It is clear that the generators of $\mathcal{C}$ pullback along $i$ to the generators of $\mathcal{B}$. The morphisms in degree zero do not change by an application of algebraic Harthogs' Lemma\footnote{Apply algebraic Harthogs' to $\Hom(S,V)^f\oplus \wedge^k S \subset X'$, the unique extension is still invariant, as the locus of invariant points is Zariski closed.}. As $\mathcal{Z}$ is an affine stack, the generators have no higher degree morphisms. By Lemma \ref{lem:tiltingbundle}, the generators of $\mathcal{B}$ also have no higher degree morphisms. Together these show fully faithfulness and therefore $\mathcal{C} \simeq \mathcal{B}$.
\end{proof}
\begin{proposition}\label{prop:nccrequivfibre}
	The subcategory $\mathcal{C}$ is equivalent to $\mathcal{D}$.
\end{proposition}
\begin{proof}
	Again, the generators of $\mathcal{C}$ pullback along $j$ to the generators of $\mathcal{D}$. We will look at local cohomology to get $j^*$ is fully faithful restricted to $\mathcal{C}$. Conceptually this is very straightforward, we are just adding a section to $\wedge^k S^* \cong \mathcal{O}(1)$ as we are removing a linear hyperplane. As $\mathcal{O}(1)$ is greater than all the generators\footnote{Using the total order on Young diagrams.} we will not get any extra morphisms.\par
	Let $\mathcal{W} = \left[\Hom(S,V) \medoplus 0/GL(S)\right]$, then we have an exact triangle
	$$\Hom(\mathcal{E}, \mathcal{F}) \rightarrow \Hom(\mathcal{E}, j_* j^*\mathcal{F}) \rightarrow \Hom\left(\mathcal{E},\mathcal{H}_{\mathcal{W}}(\mathcal{F})\right)$$
	and as $\Hom(\mathcal{E}, j_* j^*\mathcal{F}) \cong \Hom(j^*\mathcal{E}, j^*\mathcal{F})$, we need to show $\Hom\left(\mathcal{E},\mathcal{H}_{\mathcal{W}}(\mathcal{F})\right) = 0$ for all $\mathcal{E},\mathcal{F}$ in our set of generators. We have $\mathcal{H}_{\mathcal{W}}(\mathcal{F}) = \mathcal{F} \medotimes \mathcal{H}_{\mathcal{W}}(\mathcal{O}_{X})$ and 
	$$\mathcal{H}_{\mathcal{W}}(\mathcal{O}_{X}) \cong \mathcal{O}_{X} \medotimes\Sym{\bullet} (\wedge^k S) \medotimes \wedge^k S$$
	by Lemma \ref{lem:localcoh_for_subbundles}.\footnote{The $\wedge^k S$ term provides $\mathcal{O}(1)$ with a new section.} Putting this all together, let $\alpha, \beta \in \mathcal{UP}_{n,k}$, we want $$H\left(\mathcal{Z}, (\BBS^{\alpha}S^* \medotimes (\BBS^{\beta} S^*)^*) \medotimes\mathcal{H}_{\mathcal{W}}(\mathcal{O}_{\mathcal{W}})\right) = 0.$$
	So we need to show
	$$ \Sym{\bullet} \left(S \medotimes V^*\right)\medotimes \left( \medoplus_{i \geq 1} \BBS^{(i,i, \dots,i)}S\right) \medotimes\BBS^{\alpha}S^* \medotimes \left(\BBS^{\beta} S^*\right)^* $$
	has no $GL(S)$ invariant elements.\par
	Let $\BBS^{\gamma}S$ be any representation appearing in the decomposition above. We have $\gamma$ comes from some representation appearing in $(a_1,a_2, \dots , a_k) \otimes (i, \dots , i) \otimes (-\alpha_k, \dots , -\alpha_{1}) \otimes (\beta_{1}, \dots , \beta_k) $. Apply the lower bound from Lemma \ref{lem:littlewoodrichardsonestimates} to get $\gamma_{1} \geq a_k + i -\alpha_{k} + \beta_k  \geq i \geq 1$. As $\alpha_{k} = \beta_{k} =0$ by assumptions on $\alpha,\beta$ and $a_k \geq 0$.\par
	Therefore $\gamma \not = (0,0, \dots , 0)$ and there are no $GL(S)$ invariant terms.
\end{proof}
\begin{theorem}\label{thm:allNCCRsareequiv}
	Let $X = \Hom (S, V)$ and $X' = \Hom (V, Q)$. Let $\Lambda$ and $\Lambda'$ be the NCCR's of $X/SL(S)$ and $X'/SL(Q)$ respectively from Theorem \ref{thm:FGDandNCCR}. Then $\Lambda$ and $\Lambda'$ are derived equivalent.
\end{theorem}
\begin{proof}
	We have already done all the work;\par
	Propositions \ref{prop:linebundleequivfibre} and \ref{prop:nccrequivfibre} together tell us that $\Lambda$ and $\Lambda'$ both have derived categories equivalent to the respective categorical crepant resolutions, $\mathcal{B}$ and $\mathcal{B}'$. Lemma \ref{lem:differentgeneratingsets} shows that $\mathcal{B}$ and $\mathcal{B}'$ are equivalent (actually identical, just two different generating sets). Using this $\Lambda$ and $\Lambda'$ are derived equivalent as $D^b(\Lambda) \simeq \mathcal{B} = \mathcal{B}' \simeq D^b(\Lambda')$. Recall that $\mathcal{D} \simeq D^b(\Lambda)$, basically by definition.
\end{proof}

\section{Homological Projective Duality}\label{sec:HPD}
We finally want to address the question of the Homological Projective Dual to the Pl{\"u}cker embedding of the Grassmannian. This will build on the work from the previous section. For more details on HPD see the introduction, Section \ref{sec:HPDintro}.\par
We split this section into three parts. The first one deals with what happens on a single fibre. This requires us to extend Propositions \ref{prop:linebundleequivfibre} and \ref{prop:nccrequivfibre} to versions where we replace the subcategories with categories of matrix factorizations. The second piece is a summary of the results from \cite{2017arXiv170501437V} and what that gives us in our case. Finally we work globally and construct the HPD. 
\subsection{HPD for a single fibre}\label{sec:HPDfibre}
We consider the Pl{\" u}cker embedding of the Grassmannian. As before let $S$ and $V$ be vector spaces of dimensions $k$ and $n$ with $1 < k< n - 1$, assume $n$ and $k$ are coprime and also pick a trivialization of $\det V$. Let $Gr = Gr(S,V)$, we have
\[\begin{tikzcd}
{Gr} & {\BBP(\wedge^k V)}
\arrow[from=1-1, to=1-2, hook]
\end{tikzcd}\]  
and the Kuznetsov - Fonarev exceptional collection
$$D^b(Gr) = \langle \mathcal{A}, \dots ,\mathcal{A}(n-1) \rangle$$
where $\mathcal{A}$ is the set of sheaves associated to representations from $\mathcal{UP}_{n,k}$ as before and $\mathcal{O}_{Gr}(1)$ is the pullback of $\mathcal{O}_{\BBP(\wedge^k V)}(1)$. Before we approach HPD in general, we recall the motivational question.\par
In general, for $Gr_H$, which is $Gr$ intersected with a hyperplane section, we have
$$D^b(Gr_H) = \langle \mathcal{K}_H, \mathcal{A}(1), \dots, \mathcal{A}(n-1) \rangle$$
and the basic idea of HPD is to find an alternative description of the Kuznetsov component, $\mathcal{K}_H$, the "interesting" piece of the semi-orthogonal decomposition. The first step for us is to use matrix factorizations to find an alternative description of $D^b(Gr_H)$. Recall, $Gr_H$ is cut out by a section of $\mathcal{O}_{Gr}(1)$, and we have the following diagram.
\begin{center}
	\[\begin{tikzcd}
		{\pi^{-1}(Gr_H)} & {\mathcal{O}_{Gr}(-1)} \\
		{Gr_H} & Gr
		\arrow["\pi", shift left=1, from=1-2, to=2-2]
		\arrow["\pi", from=1-1, to=2-1]
		\arrow["i", hook, from=1-1, to=1-2]
		\arrow["j", from=2-1, to=2-2]
		\arrow["{\iota_0}", shift left=1, hook, from=2-2, to=1-2]
	\end{tikzcd}\]
\end{center}
Kn{\" o}rrer periodicity tells us that 
$$D^b(Gr_H) \simeq D^b(\mathcal{O}_{Gr}(-1), W)$$
where $W$ is the superpotential associated to the section $w$ cutting out $Gr_H$. As we are interested in $\mathcal{K}_H \subset D^b(Gr_H)$ we want to find a description of the equivalent subcategory of $D^b(\mathcal{O}_{Gr}(-1), W)$.\par
To do this we note that the Kn{\" o}rrer periodicity equivalence is given by $i_* \pi^*$ and as $\mathcal{K}_H$ is the right orthogonal to $\langle \mathcal{A}(1), \dots, \mathcal{A}(n-1) \rangle$ we can equivalently describe it as 
$$\mathcal{K}_H = \{\mathcal{E} | j_* \mathcal{E} \in \mathcal{A}\}.$$
Let $\iota_0 : Gr \rightarrow \mathcal{O}_{Gr}(-1)$ be the inclusion of the zero section, then we have $\iota_0^* i_* \pi^* \mathcal{E} \simeq \iota_0^* \pi^* j_* \mathcal{E} \simeq j_* \mathcal{E}$ as $\pi$ is flat and $\pi \circ \iota_0 = id$. Therefore
$$\mathcal{K}_H \simeq \mathcal{B}^W = \{\mathcal{F} | \iota_0^* \mathcal{F} \in \mathcal{A}\} \subset D^b(\mathcal{O}_{Gr}(-1), W).$$
And if we forget about $W$, we recognize $\{\mathcal{F} | \iota_0^* \mathcal{F} \in \mathcal{A}\} \subset D^b(\mathcal{O}_{Gr}(-1))$ as the subcategory $\mathcal{B}$ from Section \ref{sec:equivNCCR}.\par
As $\mathcal{B}$ is derived equivalent to $\Lambda$, the NCCR of the affine cone of the Grassmannian, we suspect that $\mathcal{K}_H \simeq D^b(\Lambda, W)$. To prove this, we need to extend $\mathcal{B} \simeq \mathcal{D}$ to $\mathcal{B}^W \simeq \mathcal{D}^W$, where $\mathcal{D}^W = D^b(A, W) \subset D^b(\mathcal{X}, W)$ and $A$ is the set of vector bundles associated to the irreducible representations of $SL(S)$ indexed by $\mathcal{UP}_{n,k}$. We will do this next.\par
First note that we can define $W$ on $\mathcal{Z} = \left[\Hom(S,V)\times \wedge^k S  /GL(S) \right]$, pick any section of $\mathcal{O}(1)$, i.e. an element $w \in \wedge^k V^*$. Recall that we are calling the coordinates of $\mathcal{Z}$, $(x,p)$. We then have $W: \mathcal{Z} \rightarrow \BBC$ defined by $W(x,p) = w(\wedge^k x \cdot p)$. The restriction of $W$ to $\mathcal{O}_{Gr}(-1)$ is exactly the same $W$ as the one used for the Kn{\" o}rrer periodicity result. We also define an $R-$charge on $\mathcal{Z}$ by letting $\BBC^*$ act with weights $0,2$ on $x,p$. This $R-$charge restricted to $\mathcal{O}_{Gr}(-1)$ is again the same as the Kn{\" o}rrer periodicity one. It also gives us an $R-$charge on $\mathcal{X} = \left[\Hom(S,V)/SL(S) \right]$.\par
Morally or conceptually, we already have $\mathcal{B} \simeq \mathcal{C} \simeq \mathcal{D}$, and these equivalences are induced by the inclusion maps $i,j$. So as we can define $W$ on $Y$ and $\mathcal{X}$ by pullback along these maps, we want to just say, add in $W$ and get $\mathcal{B}^W \simeq \mathcal{C}^W \simeq \mathcal{D}^W$ where $\mathcal{C}^W =  D^b(A, W) \subset D^b(\mathcal{Z}, W)$. In reality we need to do a little bit more, starting with alternative descriptions of these subcategories.
\begin{lemma}\label{lem:mfonstackx}
	Any object in the subcategory $\{\mathcal{E} |\mathcal{E}|_0 \in \mathcal{UP}_{n,k} \} \subset D^b(\mathcal{X}, W)$ is an element of $\mathcal{D}^W$.
\end{lemma}
\begin{proof}
	This follows directly from Lemma \ref{lem:equivalentdescriptionofwindowsubcategory}.
\end{proof}
\noindent
Exactly the same result shows that we can equivalently define $\mathcal{C}^W$ as
$$\mathcal{C}^W = \{ \mathcal{E} |\mathcal{E}|_0 \in \mathcal{UP}_{n,k} \} \subset D^b(\mathcal{Z}, W).$$
It turns out that the question of extending an equivalence of categories to an equivalence of matrix factorization categories has already been partially answered, see \cite[Section 4.2]{2014arXiv1401.3661A}. The results found there are strong enough to apply in our situation, and we will briefly explain them the first time we use them. We will also use similar ideas to prove the final HPD statement in Section \ref{sec:tautologicalHPD}.

\begin{proposition}\label{prop:nccrwithmfeqquiv}
	The subcategory $\mathcal{C}^W$ is equivalent to $\mathcal{D}^W$ by restriction along $i$.
\end{proposition}
\begin{proof}
	First, fully faithfulness follows from fully faithfulness for the case where $W=0$. This is as morphisms on $(\mathcal{Z}, W)$ can be computed using a spectral sequence whose first page is morphisms on $(\mathcal{Z}, 0)$, similarly for $\mathcal{X}$. To see this note that $\Hom(\mathcal{E}, \mathcal{F})$ is a bi-complex and therefore we get a spectral sequence whose first page is $(H^{\bullet}(\mathcal{H}om(\mathcal{E},\mathcal{F})), d_{\mathcal{E},\mathcal{F}})$ which converges to $\Hom(\mathcal{E}, \mathcal{F})$ on $(\mathcal{Z}, W)$, see \cite[Remark 2.14]{2011CMaPh.304..411S} for the details. Therefore Proposition \ref{prop:nccrequivfibre} tells us that that we have an isomorphism of spectral sequences and therefore $i^*$ is fully faithful when restricted to $\mathcal{C}^W$.\par
	In general to prove essential surjectivity, we will show that any object lifts by finding a standard form. In this case, this follows by definition of $\mathcal{D}^W$.	Let $\mathcal{F} \in \mathcal{D}^W$, by definition we can find an quasi-isomorphic matrix factorization $(F, d_F)$ where $F$ is built from vector bundles which lift to $\mathcal{C}$. As we already have fully faithfulness we get a matrix factorization, $(\tilde{F}, d_{\tilde{F}}) \in \mathcal{C}^W$, which is sent to $(F, d_F)$ by $i^*$ and therefore we have essential surjectivity.
\end{proof}
We also need an analogue of Lemma \ref{lem:mfonstackx} for $\mathcal{B}^W$. We can not restrict to the origin in $Gr$ as it does not exist, so we have to phrase things slightly differently. Let $\mathcal{H}^W  = D^b(A, W) \subset D^b(\mathcal{O}_{Gr(S,V)}(-1), W)$ be the subcategory of matrix factorizations quasi-isomorphic to those built from the vector bundles $\BBS^{\alpha} S^*$ for $\alpha \in \mathcal{UP}_{n,k}$.
\begin{proposition}\label{prop:linebundlevectorbundlegen}
	The subcategories $\mathcal{B}^W$ and $\mathcal{H}^W$ are equal.
\end{proposition}
\begin{proof}
	First note that $\mathcal{H}^W \subset \mathcal{B}^W$ as any matrix factorization built from the $\BBS^{\alpha} S^*$ satisfies the restriction condition. Also $\mathcal{H}^W$ is an admissible subcategory as $\End(\medoplus \BBS^{\alpha}S^*) \cong \Lambda$ by Theorem \ref{thm:allNCCRsareequiv} and $\Lambda$ has finite global dimension by Theorem \ref{thm:FGDandNCCR}.\par
	Now let $\mathcal{F} \in D^b(\mathcal{O}_{Gr}(-1), W)$, we can restrict  $\mathcal{F}$ to $Gr$ and then project into $\mathcal{A}$ or as these commute, we can first project into $\mathcal{H}^W$ and then restrict to $Gr$. This is as the projection functor to $\mathcal{A}$ is the restriction of the projection functor to $\mathcal{H}^W$. In particular let $\mathcal{F} \in \mathcal{B}^W$ and then let $\mathcal{F}^{\perp}$ be projection into the orthogonal to $\mathcal{H}^W$. By definition of $\mathcal{B}^W$ we have that $i^*\mathcal{F}^{\perp} \in \mathcal{A}$. We can also first restrict to $Gr$ and then project. By definition we get 0 and therefore $i^*\mathcal{F}^{\perp} \simeq 0$.	As $\mathcal{F}^{\perp}$ is equivariant this tells us that in fact $\mathcal{F}^{\perp} \simeq 0$ and $\mathcal{B}^W = \mathcal{H}^W$. 
\end{proof}
\begin{theorem}\label{thm:linebundleenccrquivpotential}
	$\mathcal{B}^W \simeq \mathcal{D}^W$
\end{theorem}
\begin{proof}
	We have already seen that $\mathcal{D}^W \simeq \mathcal{C}^W$, so we just have to show that $\mathcal{B}^W \simeq \mathcal{C}^W$ and we have already done all the required work to show that.\par
	Fully faithfulness follows by exactly the same argument as in Proposition \ref{prop:nccrwithmfeqquiv}, this time using Proposition \ref{prop:linebundleequivfibre} for the $W=0$ case.	Then Proposition \ref{prop:linebundlevectorbundlegen} shows that $j^*:\mathcal{C}^W \rightarrow \mathcal{B}^W$ is essentially surjective, as essentially by definition, every object of $\mathcal{H}^W = \mathcal{B}^W$ lifts to an object of $\mathcal{C}^W$.
\end{proof}
While in some sense this result is the heart of the HPD statement, it is not the full result. For that we need to define a subcategory over $\BBP(\wedge^k V^*)$ that when restricted to each fibre is equivalent to $\mathcal{D}^W$. Then we need to prove all the equivalences of subcategories for linear subspaces, $L \subset \wedge^k V^*$.\par
It turns out that this splits into two parts. The first follows from a general theorem of Rennemo's, \cite[Thm 1.1]{2017arXiv170501437V}, this part can be seen as a generalization of replacing $\mathcal{K}_H$ with the matrix factorization category $\mathcal{B}^W$. We will sketch this argument in the next section. The second part is proving a generalization of Theorem \ref{thm:linebundleenccrquivpotential} which will be done in the section after.
\subsection{Tautological HPD}\label{sec:tautologicalHPD}
This section is all the work of Rennemo ,\cite{2017arXiv170501437V}, he works in a greater generality, we quote his results as they apply in our situation. \par
As in the introduction, let $L$ be a linear subspace of $\wedge^k V^* = H^0(Gr, \mathcal{O}(1))$ and $L^{\perp} \subset \wedge^k V$ be its annihilator. We have a natural superpotential, $\mathcal{O}_{Gr} (1) \medotimes L \rightarrow \BBC$ given by evaluation and the standard $R-$charge for vector bundles. 
\begin{proposition}[{\cite[Prop. 2.6]{2017arXiv170501437V}}]\label{prop:matrixequiv}
	Let $W$ be the natural function $\mathcal{O}(-1)\medotimes L \rightarrow \BBC$, and let $w$ be the induced section on $\mathcal{O}(1)\medotimes L$. Then $D^b(\mathcal{O}(-1)\medotimes L, W) \cong D^b(Gr_{L^{\perp}})$ as long as $Gr_{L^{\perp}}$ has the expected dimension, where $Gr_{L^{\perp}} = \{w=0\}$ is the intersection of $Gr$ with $L^{\perp}$ inside $\BBP(\wedge^k V)$.
\end{proposition}
\noindent
We rewrite in terms of global quotient stacks. Let $\Hom (S,V)^f \subset \Hom(S,V)$ be the set of full rank maps, then $Gr \cong \left[ \Hom(S,V)^f/GL(S) \right]$ and $\mathcal{O}_{Gr}(-1) \cong \left[ \Hom(S,V)^f \times \wedge^k S/GL(S) \right]$. Rennemo then realises $\mathcal{O}_{Gr}(-1)\medotimes L$ as a solution to a GIT problem. Consider \\$\mathcal{Z}' = \left[\Hom(S,V)^f \times \wedge^k S \times L/GL(S)\times \BBC^*\right]$ where $\BBC^*$ acts with weights $0,1,-1$. One stability condition has unstable locus $\Hom(S,V)^f \times 0 \times L$ and the GIT quotient is isomorphic to the total space of $\mathcal{O}_{Gr}(-1)\medotimes L$. The other stability condition has unstable locus $\Hom(S,V)^f \times \wedge^k S \times 0$ and the GIT quotient is  $Y_L \cong \mathcal{O}_{Gr\times \BBP L}(-1,-1)$.\par
Note that $\mathcal{Z}'$ still has a natural potential, $W$, which is effectively evaluation. If we label the variables by $(x,p,l)$ then $W(x,p,l) = l(\wedge^k x (p))$. This potential restricts to both GIT quotients, on $\mathcal{O}_{Gr}(-1)\medotimes L$ it is exactly the potential from earlier proposition. On the other side, $Y_L$, the associated section cuts out the incidence hyperplane $\{(x,l) | l(x) = 0\} \subset Gr \times \BBP L$. We also define an $R-$charge which acts with weights $0,2,0$. Finally we set 
$$\mathcal{B}^W_L = \{\mathcal{E} | \mathcal{E}|_{Gr\times [l]} \in \mathcal{A}, \forall 0 \not=l \in L\} \subset D^b(Y_L, W)$$
where $Gr \times [l]$ is included as the zero section and 
$$\mathcal{D}_{L^{\perp}} = D^b(\mathcal{O}_{Gr}(-1)\medotimes L, W).$$
By Proposition \ref{prop:matrixequiv}, $\mathcal{D}_{L^{\perp}}$ is the replacement for $D^b(Gr_{L^{\perp}})$. The following is then the main result from \cite{2017arXiv170501437V}.
\begin{theorem}[{\cite[Thm 1.1]{2017arXiv170501437V}}]\label{thm:firstequiv}
	$\mathcal{B}^W_{\wedge^k V^*}$ is an homological projective dual to $Gr(S,V) \hookrightarrow \BBP (\wedge^k V) $.
\end{theorem}
Recall, $\mathcal{B}^W_{\wedge^k V^*}$ is a category linear over $\BBP (\wedge^k V^*)$. Ideally we would wish it to be equivalent to $D^b(Y')$ for some variety $Y' \rightarrow \BBP (\wedge^k V^*)$, but as discussed in Section \ref{sec:HPDintro}, this is not the case.
\begin{remark}
	Considering the case $\dim L = 1$, we get $Y_L = \mathcal{O}_{Gr}(-1)$ and by the argument of the previous section, $\mathcal{B}^W_L \simeq \mathcal{K}_H$. This provides motivation for the definition of $\mathcal{B}^W_L$.
\end{remark}

\subsection{Geometric HPD}\label{sec:geometricHPD}
One can interpret $\mathcal{B}^W_L$ as being the tautological or canonical HP dual, it always exists in large generality. However in our situation, described the way it currently is, it does not look very geometric. To get a more geometric HP dual, we will use another GIT problem, this time considering the group $GL(S)$, this is a generalization of what we did in Section \ref{sec:equivNCCR}. We will use this to show another equivalence of categories which extends the earlier result for fibres. All of the following stacks are over $\BBP L$.\par
We have
$$ \mathcal{Z}_L = \left[\Hom(S,V)\times \wedge^k S \times (L\setminus 0)/GL(S)\times \BBC^*_{0,-1,1} \right].$$
Call the coordinates $(x,p,l)$. One GIT quotient is 
$$Y_L = \left[\Hom(S,V)^f\times \wedge^k S \times (L\setminus 0)/GL(S)\times \BBC^*_{0,-1,1} \right]$$
as above, giving us $\mathcal{O}_{Gr\times \BBP L}(-1,-1)$, the space we use to define $\mathcal{B}^W_L$. The other quotient is 
$$\mathcal{X}_L = \left[\Hom(S,V)\times (\wedge^k S\setminus 0) \times (L\setminus 0)/GL(S)\times \BBC^*_{0,-1,1} \right].$$
This stack also maps down to $\BBP L$ and the fibres are isomorphic to $\mathcal{X} = [\Hom(S,V)/SL(S)]$. 
\begin{center}
	\[\begin{tikzcd}
	& {\mathcal{Z}_L =\left[\frac{\Hom(S,V)\times \wedge^k S \times (L\backslash 0)}{GL(S)\times \BBC^*_{0,-1,1}} \right]} \\
	{Y_L =\left[\frac{\Hom(S,V)^f\times \wedge^k S \times (L\backslash 0)}{GL(S)\times \BBC^*_{0,-1,1}} \right]} & {} & {\mathcal{X}_L = \left[\frac{\Hom(S,V)\times (\wedge^k S\backslash 0) \times (L\backslash 0)}{GL(S)\times \BBC^*_{0,-1,1}} \right]} \\
	{\mathcal{O}_{Gr\times \BBP L}(-1,-1)} & {\BBP L} & {\mathcal{X} \hookrightarrow \mathcal{X}_L \twoheadrightarrow\BBP L}
	\arrow["\cong"{description}, from=3-1, to=2-1]
	\arrow["\cong"{description}, from=3-3, to=2-3]
	\arrow["i", hook, from=2-1, to=1-2]
	\arrow["j"', hook', from=2-3, to=1-2]
	\arrow[from=1-2, to=3-2]
	\arrow[from=2-1, to=3-2]
	\arrow[from=2-3, to=3-2]
	\end{tikzcd}\]
\end{center}

To extend Theorem \ref{thm:linebundleenccrquivpotential} to the general case over $\BBP(L)$ we will use a similar argument. \par
First we need to define the more general subcategories which we will do with a type of grade restriction rule, like the definition of $\mathcal{B}^W_L$. To prove fully faithfulness it will be sufficient to show it on fibres, and we have already seen that result. For essential surjectivity we will again prove that any object of the subcategory has a representation that obviously lifts to $\mathcal{Z}_L$.\par
Consider the $\mathcal{X}_L$ side, on a fibre we have the stack $\mathcal{X}$ which is where the NCCR, $\Lambda$, is found. More formally, let $0 \not =l \in L$, we have the line $[l]$ and an inclusion map,
$$ \mathcal{X} \cong \left[\frac{\Hom(S,V) \times \wedge^k S \setminus 0}{GL(S)}\right] \cong \left[\frac{\Hom(S,V) \times \wedge^k S\setminus 0\times [l] \setminus 0}{GL(S)\times \BBC^*}\right] \hookrightarrow \mathcal{X}_L.$$
Call this map $\iota_{[l]}$, then we define
$\mathcal{D}_L^W\subset D^b(\mathcal{X}_L, W)$ as follows
$$\mathcal{D}_L^W = \{\mathcal{E} | \iota_{[l]}^* \mathcal{E} \in \mathcal{D}^W \;\;\; \forall 0 \not = l \in L  \}. $$
We can describe this subcategory equivalently using Lemma \ref{lem:mfonstackx}, which tells us that $\mathcal{D}^W \subset D^b(\mathcal{X}, W)$ can be equivalently described as the matrix factorizations $\mathcal{E}$ such that $\mathcal{E}|_0 \in \mathcal{UP}_{n,k}$ and therefore we can say
$$\mathcal{D}_L^W = \{\mathcal{E} | \mathcal{E}_{0\times[l]} \in \mathcal{UP}_{n,k} \;\;\; \forall 0 \not = l \in L  \}.$$
This definition is a natural way to extend $D^W$ over $\BBP L$, but it is not a convenient description for proving essential surjectivity. For that we define the subcategory $\mathcal{H}_{\mathcal{X}_L} \subset D^b(\mathcal{X}_L, W)$ as all matrix factorizations homotopic to those built from $\BBS^{\alpha}S^* \medotimes \pi^* \mathcal{E}$ where $\alpha \in \mathcal{UP}_{n,k}$ and $\mathcal{E}$ is any sheaf on $\BBP L$. (One can think of $\BBS^{\alpha}S^*$ as coming from $\left[\Hom(S,V)/GL(S) \right]$.) This is a $\BBP L$-linear category, and clearly any object in $\mathcal{H}_{\mathcal{X}_L}$ lifts to an object in $D^b(\mathcal{Z}_L, W)$.
\begin{proposition}\label{prop:admissible}
	The subcategory $\mathcal{H}_{\mathcal{X}_L}$ is right admissible.
\end{proposition}
\begin{proof}
	Let $\mathcal{F} = \medoplus_{\alpha \in \mathcal{UP}_{n,k}} \BBS^{\alpha} S^* \medotimes \mathcal{O}$. Then the projection functor is given by $\mathcal{E} \mapsto  \mathcal{H}om_{\BBP L}(F, \mathcal{E})$ with adjoint $\mathcal{F} \medotimes_{\mathcal{E} nd (F)} -$. (Note that we could modify each summand of $\mathcal{F}$ by any vector bundle from $\BBP L$, the endomorphism algebra would then be Morita equivalent.) For this to all work we only need to check that $\mathcal{E} nd (F)$ has finite global dimension. This follows as on an affine open, $U$, we have $\mathcal{E} nd (F)(U) \cong \Lambda \otimes R$ and both $\Lambda$ and $R$ have finite global dimension. (For example, if we take the standard affine charts on $\BBP L$, $R = \BBC[x_1, \dots x_{\dim L -1}]$.)
\end{proof}
\begin{proposition}\label{prop:xltwodescriptions}
	The subcategories $\mathcal{D}_L^W$ and $\mathcal{H}_{\mathcal{X}_L}$ are equal.
\end{proposition}
\begin{proof}
	Clearly $\mathcal{H}_{\mathcal{X}_L} \subset \mathcal{D}_L^W$ so we just need to prove the other inclusion. 	$\mathcal{H}_{\mathcal{X}_L}$ is right admissible by Proposition \ref{prop:admissible}, so we can project $\mathcal{F} \in \mathcal{D}_L^W$ into the right orthogonal of $\mathcal{H}_{\mathcal{X}_L}$. It is then sufficient to prove this projection is zero.\par
	First restricting to any fibre, and then projecting, it is clear that this projection is zero, by Lemma \ref{lem:equivalentdescriptionofwindowsubcategory}. Restricting and projecting commute, so therefore the projection is zero on each fibre, and therefore zero everywhere.
\end{proof}
Next we look at the $Y_L$ side, it is essentially the same as the $\mathcal{X}_L$ side. We have
$$\mathcal{B}_L^W = \{\mathcal{E} | \mathcal{E}|_{Gr\times [l]} \in \mathcal{A}, \forall 0 \not=l \in L\} \subset D^b(Y_L, W)$$
and Proposition \ref{prop:linebundlevectorbundlegen} gives us the equivalent definition,
$$\mathcal{B}_L^W = \{\mathcal{E} | \mathcal{E}|_{\mathcal{O}_{Gr}(-1)\times [l]} \in \mathcal{B}^W, \forall 0 \not=l \in L\} \subset D^b(Y_L, W).$$
To prove essential surjectivity, again we would like a different description. Let $\mathcal{H}_{Y_L} \subset D^b(\mathcal{Y}_L, W)$ be defined as all matrix factorizations built from $\BBS^{\alpha}S^* \medotimes \pi^* \mathcal{E}$ where $\alpha \in \mathcal{UP}_{n,k}$ and $\mathcal{E}$ is any sheaf on $\BBP L$. (One can think of $\BBS^{\alpha}S^*$ as coming from $Gr$.) This is a $\BBP L$ linear category, and clearly any object in $\mathcal{H}_{Y_L}$ lifts to an object in $D^b(\mathcal{Z}_L, W)$.
\begin{proposition}\label{prop:yltwodescriptions}
		The subcategories $\mathcal{B}^W_{L}$ and $\mathcal{H}_{Y_L}$ are equal.
\end{proposition}
\begin{proof}
	Identical to Proposition \ref{prop:xltwodescriptions}.
\end{proof}

The final preparation step is to define $\mathcal{C}_L^W \subset D^b(\mathcal{Z}_L, W)$ as all matrix factorizations built from the sheaves $\BBS^{\alpha}S^* \medotimes \pi^* \mathcal{E}$ as before.\footnote{We could instead define $\mathcal{C}^W_L$ as matrix factorizations $\mathcal{E}$ such that $\mathcal{E}|_{0 \times [l]}$ is composed of representations from $\mathcal{UP}_{n,k}$. An analogue of Proposition \ref{prop:xltwodescriptions} would prove these definitions are equivalent.} 
Then we can finally prove the final steps of the HPD.\\
\begin{theorem}\label{thm:secondequiv}
	The subcategories $\mathcal{B}_L^W$ and $\mathcal{D}_L^W$ are equivalent.
\end{theorem}
\begin{proof}
	We will prove this in a very similar manner to Theorem \ref{thm:linebundleenccrquivpotential}, by proving that both subcategories are equivalent to $\mathcal{C}_L^W$, via the maps $i^*$ and $j^*$.\par
	Fully faithfulness is straightforward. Let $\mathcal{E}, \mathcal{F} \in \mathcal{C}_L^W$, then consider morphisms between them as a sheaf on $\BBP L$, $\mathcal{H}om_{\BBP L} (\mathcal{E}, \mathcal{F})$. The functor $i^*$ then induces a map of sheaves
	$$\mathcal{H}om_{\BBP L} (\mathcal{E}, \mathcal{F}) \rightarrow \mathcal{H}om_{\BBP L} (i^*\mathcal{E}, i^*\mathcal{F}).$$
	We want that this map is an isomorphism and so it is sufficient to prove fully faithfulness on each fibre.  Proposition \ref{prop:nccrwithmfeqquiv} tells us this. The same argument, this time using Theorem \ref{thm:linebundleenccrquivpotential} for the result on fibres, gives the result for $j^*$.\par
	It is clear that the images of $\mathcal{C}_L^W$ under $j^*$ and $i^*$ are $\mathcal{H}_{\mathcal{X}_L}$ and $\mathcal{H}_{Y_L}$ respectively, it is also clear that every element in these subcategories are mapped to. Propositions \ref{prop:xltwodescriptions} and \ref{prop:yltwodescriptions} tell us that these subcategories are actually $\mathcal{D}_L^W$ and $\mathcal{B}_L^W$.\par
	Putting this all together gives us $\mathcal{B}_L^W \simeq \mathcal{C}^W_L \simeq \mathcal{D}_L^W$.
\end{proof}
\begin{remark}
	This theorem would be a special case of a more general conjecture. Let $X \rightarrow Y$ and let $\mathcal{C} \subset D^b(X)$, $\mathcal{D} \subset D^b(Y)$ such that $\mathcal{C} \simeq \mathcal{D}$. Then let $Z$ be any smooth scheme and let $X_Z, Y_Z$ be schemes over $Z$ such that every fibre is isomorphic to $X,Y$. We can set $\mathcal{C}_Z = \{\mathcal{E}| \mathcal{E}|_{X\times z} \in \mathcal{C} \; \forall z\}$ and define $\mathcal{D}_Z$ similarly. Then is $\mathcal{C}_Z \simeq \mathcal{D}_Z$?
\end{remark}
We now have everything we need and can define the HPD, let $\mathcal{D}^W_{\Lambda} = \mathcal{D}^W_{\wedge^k V^*}$.\\
\begin{theorem}\label{thm:finalequiv}
	$\mathcal{D}^W_{\Lambda}$ is a homological projective dual for $Gr(S,V) \rightarrow \BBP(\wedge^k V)$.
\end{theorem}
\begin{proof}
	Combine Theorem \ref{thm:firstequiv} and Theorem \ref{thm:secondequiv} to get the result.
\end{proof}
In words, the homological projective dual of the Pl{\" u}cker embedding of the Grassmannian is the global version of the NCCR for the affine cone of the Grassmannian with a natural superpotential.
\begin{remark}
	In principle we should not need to prove Theorem \ref{thm:finalequiv} by splitting into two different GIT problems. We could instead work on the larger stack \\ $\left[\Hom(S,V) \times \wedge^k S \times L/GL(S) \times \BBC^*\right]$ and define various window subcategories here. This would however require us to redo most of \cite{2017arXiv170501437V}. So instead we followed the path seen above, pulling on the full power of Rennemo's results. One can think of Theorem \ref{thm:firstequiv} as containing the technical HPD results and then Theorem \ref{thm:secondequiv} provides the geometrical interpretation.
\end{remark}

\newpage

\printbibliography
\end{document}